\documentclass[11pt]{article}

\usepackage{xypic}
\usepackage{amsgen}
\usepackage{amsmath}
\usepackage{amstext}
\usepackage{amsbsy}
\usepackage{amsopn}
\usepackage{amsfonts}
\usepackage{amssymb}
\usepackage{graphicx}
\usepackage{pstricks}
\xyoption{all}

\def\Box{\square}

\def\mapright#1{\smash{\mathop{\longrightarrow}\limits^{#1}}}
\def\tra#1{\smash{\mathop{\mid\kern
-1pt\joinrel\relbar\joinrel\relbar}\limits^{*}_{#1}}}
\def\longtra#1{\smash{\mathop{\mid\kern
-1pt\joinrel\relbar\joinrel\relbar\joinrel\relbar}\limits^{*}_{#1}}}
\def\vlongtra#1{\smash{\mathop{\mid\kern
-1pt\joinrel\relbar\joinrel\relbar\joinrel\relbar\joinrel\relbar}\limits^{*}_{#1}}}
\def\vvlongtra#1{\smash{\mathop{\mid\kern
-1pt\joinrel\relbar\joinrel\relbar\joinrel\relbar\joinrel\relbar\joinrel\relbar}\limits^{*}_{#1}}}
\def\vvvlongtra#1{\smash{\mathop{\mid\kern
-1pt\joinrel\relbar\joinrel\relbar\joinrel\relbar\joinrel\relbar\joinrel\relbar\joinrel\relbar}\limits^{*}_{#1}}}
\def\etra#1{\smash{\mathop{\mid\kern
-1pt\joinrel\relbar\joinrel\relbar}\limits_{#1}}}
\def\mapleft#1{\smash{\mathop{\longleftarrow}\limits^{#1}}}

\def\vlongrightarrow{\relbar\joinrel\longrightarrow}
\def\vvlongrightarrow{\relbar\joinrel\vlongrightarrow}

\def\vlongmapright#1{\smash{\mathop{\vvlongrightarrow}\limits^{#1}}}

\def\A{{\cal{A}}}
\def\iff{\Leftrightarrow}
\def\Rw{\Rightarrow}
\def\oo{\overline}

\def\wt{\widetilde}
\def\B{{\cal{B}}}
\def\C{{\cal{C}}}

\newcommand{\N}{{\rm I}\kern-2pt {\rm N}}

\def\S{{\cal{S}}}
\def\rat{\mbox{Rat}\,}

\def\pref{\mbox{Pref}}

\def\D{{\cal{D}}}

\def\p{\varphi}

\def\inv{^{-1}}

\def\bi{\begin{itemize}}
\def\ei{\end{itemize}}
\def\beq{\begin{equation}}
\def\eeq{\end{equation}}

\xyoption{all}

\newtheorem{T}{Theorem}[section]
\newcommand{\bt}{\begin{T}}
\newcommand{\et}{\end{T}}
\newcommand{\ftd}{$\square$\end{T}}

\newtheorem{Proposition}[T]{Proposition}
\newcommand{\bp}{\begin{Proposition}}
\newcommand{\ep}{\end{Proposition}}
\newcommand{\fpd}{$\square$\end{Proposition}}

\newtheorem{Lemma}[T]{Lemma}
\newcommand{\bl}{\begin{Lemma}}
\newcommand{\el}{\end{Lemma}}
\newcommand{\fld}{$\square$\end{Lemma}}

\newtheorem{Corol}[T]{Corollary}
\newcommand{\bc}{\begin{Corol}}
\newcommand{\ec}{\end{Corol}}
\newcommand{\fcd}{$\square$\end{Corol}}

\newtheorem{Result}[T]{Result}
\newcommand{\br}{\begin{Result}}
\newcommand{\er}{\end{Result}}
\newcommand{\frd}{$\square$\end{Result}}

\newtheorem{Example}[T]{Example}
\newcommand{\be}{\begin{Example}}
\newcommand{\ee}{\end{Example}}

\newtheorem{Problem}[T]{Problem}
\newcommand{\bq}{\begin{Problem}}
\newcommand{\eq}{\end{Problem}}

\newcommand{\proof}
   {\par\medbreak\noindent{\bf Proof}.\enspace}

\newcommand{\qed}{
$\Box$
\par\bigbreak}

\normalbaselines
\def\abstract#1{\par\bigskip
\begingroup\small
\baselineskip=12truept
\begin{center}ABSTRACT\end{center}
\par\medskip\par\noindent
\null\hfill\hbox{\vbox{\hsize=5truein\noindent#1}}
\hfill\null\par\endgroup\par}

\title{Finite automata for Schreier graphs of
  virtually free groups}
\author{{\bf Pedro V. Silva}\\ $ $\\ {\em Centro de
Matem\'{a}tica, Faculdade de Ci\^{e}ncias, Universidade do
Porto,}\\ {\em R. Campo Alegre 687, 4169-007 Porto, Portugal}\\
{\em e-mail:} pvsilva@fc.up.pt\\
$ $\\
{\bf Xaro Soler-Escriv\`a}\\ $ $\\ {\em Dpt. d'Estad\'{i}stiques i Inv. Operativa,
Universitat d'Alacant,}\\ {\em Apartat de correus 99,
03080-Alacant, Spain}\\
{\em e-mail:} xaro.soler@ua.es\\
$ $\\
{\bf Enric Ventura}\\ $ $\\ {\em Dept. Mat. Apl. III, Universitat
  Polit\`ecnica de Catalunya,}\\ {\em Manresa, Barcelona, Catalunya}\\
{\em e-mail:} enric.ventura@upc.edu}
\date{\today}

\begin{document}
\maketitle

\begin{center}\small
2010 Mathematics Subject Classification: 20F10, 20E06, 20-04

\medskip

Keywords: virtually free groups, Stallings foldings, Schreier graphs,
finitely generated subgroups
\end{center}

\abstract{The Stallings construction for f.g. subgroups of free groups is generalized by introducing the concept of
Stallings section, which allows an efficient computation of the core of a Schreier graph based on edge folding. It is
proved that those groups admitting Stallings sections are precisely f.g. virtually free groups, through a constructive
approach based on Bass-Serre theory. Complexity issues and applications are also discussed.}

\section{Introduction}

Finite automata became over the years the standard representation of
finitely generated subgroups $H$ of a free group $F_A$. The {\em Stallings
  construction} constitutes a simple and efficient algorithm for
building an automaton $\S(H)$ which can be used for solving the
membership problem of $H$ in $F_A$ and many other applications. This
automaton $\S(H)$ is nothing more than the {\em core automaton} of the
Schreier graph (automaton) of $H$ in $F_A$, whose structure can be described as
$\S(H)$ with finitely many infinite trees adjoined. Many features of
$\S(H)$ were (re)discovered over the years and were known to
Reidemeister, Schreier, and particularly Serre \cite{Ser}. One of the greatest
contributions of Stallings \cite{Sta} is certainly the algorithm to construct
$\S(H)$: taking a finite set of generators $h_1,\ldots,h_m$ of $H$ in
reduced form, we start with the so-called flower automaton, where
{\em petals} labelled by the words $h_i$ (and their inverse edges) are
glued to a basepoint $q_0$:
$$
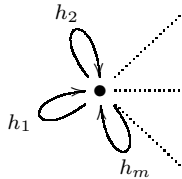
\xymatrix{
& \ \\
\bullet
\ar@(dl,l)^{h_1}
\ar@(ul,u)^{h_2}
\ar@{.}[ur] \ar@{.}[r] \ar@{.}[dr] \ar@(dr,d)^{h_m} & \ \\
& \ }$$
Then we proceed by successively folding pairs of edges of the form $q \mapleft{a} p \mapright{a} r$ until no
more folding is possible (so we get an inverse automaton). And we will have just built $\S(H)$. For details and
applications of the Stallings construction, see \cite{BS1, KM, MVW}.

Since $\S(H)$ turns out to be the core of the Schreier graph of $H\leq F_A$, this construction is independent of the
finite set of generators of $H$ chosen at the beginning, and of the particular sequence of foldings followed. And the
membership problem follows from the fact that $\S(H)$ recognizes all the reduced words representing elements of $H$...
and the reduced words constitute a section for any free group.

Such an approach invites naturally generalizations for further
classes of groups. For instance, an elegant geometric construction of
Stallings type
automata was achieved for amalgams of finite groups by Markus-Epstein
\cite{ME}. On the other hand, the most general results were obtained
by Kapovich, Weidmann and Miasnikov \cite{KWM}
for finite graphs of groups where each vertex group is either
polycyclic-by-finite or word-hyperbolic and locally quasiconvex, and
where all edge groups are virtually polycyclic. However, the complex
algorithms were designed essentially to solve the generalized word
problem, and it seems very hard to extend other features of the free
group case, either geometric or algorithmic. Our goal in the present
paper is precisely to develop a Stallings type approach with some generality
which is robust enough to exhibit several prized algorithmic and
geometric features, namely in connection with Schreier graphs.
Moreover, we succeed on identifying
those groups $G$
for which it can be carried on: (finitely generated) virtually free
groups.

Which ingredients shall we need to get a Stallings type algorithm?
First of all, we need a section $S$ with good properties that may emulate
the role played by the reduced words in the free group. In particular,
we need a rational language (i.e. recognizable by a finite automaton). We may of
course need to be more restrictive than taking all reduced words, if
we want our finite automaton to
recognize all the representatives of $H \leq_{f.g.} G$ in $S$. To get
inverse automata, it is also convenient to have $S = S\inv$

Second, the set $S_g$ of words of $S$ representing a certain $g \in G$ must be at least rational, so we can get a
finite automaton to represent each of the generalized petals.

Third, the folding process to be performed in the (generalized) flower
automaton (complemented possibly by other identification operations)
must ensure
in the end that all representatives of elements
of $H$ in $S$ are recognized by the automaton. And folding is the
automata-theoretic translation of the reduction process $w \to \oo{w}$
taking place in the free group. So we need the
condition $S_{g_1g_2} \subseteq \oo{S_{g_1}S_{g_2}}$, to make sure
that the petals (corresponding to the generators of $H$) carry enough
information to produce, after the subsequent folding, all the
representatives of elements of $H$. And this is how we were led to our
definition of {\em Stallings section}.

It is somewhat surprising how much we can get from this concept, that
turned out to be more robust than one would expect. Among other
features, we can mention independence from the generating set (so we
can have Stallings automata for free groups when we consider a non
canonical generating set!), or a generalized version of the classical
Benois Theorem. We present some applications of the whole theory,
believing that many others should follow in due time, as it happened
in the free group case.

The paper is structured as follows. In Section 2 we present the
required basic concepts. The theory of Stallings sections is presented
in Section 3. In Section 4, we discuss the complexity of the
generalized Stallings
construction in its most favourable version.
In Sections 5 and 6 we show that existence of a
Stallings section is inherited through free products with amalgamation
over finite groups
and HNN extensions over finite groups, respectively. In Section 7, we
prove that those groups admitting a Stallings section are precisely
the finitely generated virtually free groups. In Section 8, we show
that we can assume stronger properties for Stallings sections with an
eye to applications, namely the characterization of finite index
subgroups. Finally, we present some examples in Section 9.

\section{Preliminaries}

Given a finite alphabet $A$, we denote by $A^*$ the {\em free monoid
  on} $A$, with 1 denoting the empty word. A subset of a free monoid
is called a {\em language}.

We say that $\A = (Q,q_0,T,E)$ is a (finite) $A$-automaton if:
\begin{itemize}
\item
$Q$ is a (finite) set;
\item
$q_0 \in Q$ and $T \subseteq Q$;
\item
$E \subseteq Q \times A \times Q$.
\end{itemize}
A {\em nontrivial path} in $\A$ is a sequence
$$p_0 \mapright{a_1} p_1 \mapright{a_2} \,\,\cdots \,\, \mapright{a_n} p_n$$
with $(p_{i-1},a_i,p_i) \in E$ for $i = 1,\ldots,n$. Its {\em label} is the word $a_1\cdots a_n \in A^+ = A^* \setminus
\{ 1\}$. It is said to be a {\em
  successful} path if $p_0 = q_0$ and $p_n \in T$. We consider also
the {\em trivial path} $p \mapright{1} p$ for $p \in Q$. It is
successful if $p = q_0 \in T$. The {\em language} $L(\A)$ {\em
  recognized by} $\A$
is the set of all labels of successful paths in $\A$. A path of
minimal length between two vertices is called a {\em geodesic}, and so
does its label by extension.

The automaton $\A = (Q,q_0,T,E)$ is said to be {\em deterministic} if,
for all $p \in Q$ and $a \in A$, there is at most one edge of the form
$(p,a,q)$.
We say that $\A$ is {\em trim} if
every $q \in Q$ lies in some successful path.

Given deterministic
$A$-automata $\A = (Q,q_0,T,E)$ and $\A' = (Q',q'_0,T',E')$, a
morphism $\p:\A \to \A'$ is a mapping $\p: Q \to Q'$ such that
\bi
\item
$q_0\p = q'_0$ and $T\p \subseteq T'$;
\item
$(p\p, a, q\p) \in E'$ for every $(p,a,q) \in E$.
\ei
It follows that $L(\A) \subseteq L(\A')$ if there is a morphism
$\p:\A \to \A'$. The morphism $\p:\A \to \A'$ is:
\bi
\item
{\em injective} if it is
injective as a mapping $\p: Q \to Q'$;
\item
an {\em isomorphism} if it is injective, $T' = T\p$  and every edge of $E'$
is of the form $(p\p, a, q\p)$ for some $(p,a,q) \in E$.
\ei

The
{\em star} operator on $A$-languages is defined by
$$L^* = \bigcup_{n\geq 0} L^n,$$
where $L^0 = \{ 1 \}$.
A language $L \subseteq A^*$ is said to be {\em rational} if $L$ can be obtained
from finite languages using finitely many times the operators union,
product and star (admits a {\em rational expression}). Alternatively,
$L$ is rational if and only if it is
recognized by a finite (deterministic)
$A$-automaton $\A = (Q,q_0,T,E)$ \cite[Section III]{Ber}. The
definition generalizes to subsets of an arbitrary monoid in the obvious way.

We denote the set of all rational languages $L \subseteq A^*$ by $\rat
A^*$. Note that $\rat A^*$, endowed with the product of languages,
constitutes a monoid.


In the statement of a result, we shall say that a rational language $L$ is {\em
  effectively constructible} if there exists an algorithm to produce
from the data
implicit in the statement a
finite $A$-automaton $\A$ recognizing $L$.

It is convenient to summarize some closure and decidability properties of
rational languages in the following result (see \cite{Ber}
e.g.).
The prefix set of a language $L \subseteq A^*$ is defined as
$$\pref(L) = \{ u \in A^* \mid uA^* \cap L \neq \emptyset \}.$$
A {\em rational substitution} is a morphism $\p:A^* \to \rat B^*$ (where $\rat B^*$ is endowed with the product of
languages). Given $K \subseteq A^*$, we denote by $K\p$ the language $\cup_{u \in K} u\p \subseteq B^*$. Since
singletons are rational languages, monoid homomorphisms constitute particular cases of rational substitutions.

\bp
\label{proprat}
Let $A$ be a finite alphabet and let $K,L \subseteq A^*$ be rational. Then:
\bi
\item[(i)] $K\cup L, K \cap L, A^* \setminus L, \pref(L)$
  are rational;
\item[(ii)] if $\p:A^* \to \rat B^*$ is a rational substitution, then
  $K\p$ is rational; 
\item[(iii)] if $\p:A^* \to M$ is a monoid homomorphism and $M$ is
  finite, then $X\p\inv$ is rational for every $X \subseteq M$.
\ei
Moreover, all the constructions are effective, and the inclusion $K
\subseteq L$ is decidable.
\ep

Given an $A$-automaton $\A$ and $L \subseteq A^*$, we denote by $\A
\sqcap L$ the $A$-automaton obtained by removing from $\A$ all the
vertices and edges which do not lie in some successful path
labelled by a word in $L$.

\bp
\label{clsq}
Let $\A$ be a finite $A$-automaton and let $L \subseteq A^*$ be a
rational language. Then $\A
\sqcap L$ is effectively constructible.
\ep

\proof
Write $\A = (Q,q_0,T,E)$ and let $\A' = (Q',q'_0,T',E')$ be a finite
$A$-automaton recognizing $L$. The {\em direct product}
$$\A'' = (Q \times Q', (q_0,q'_0), T \times T', E'')$$ is defined by
$$E'' = \{ ((p,p'),a,(q,q')) \mid (p,a,q) \in E,\; (p',a,q') \in E'
\}.$$ Let $\B$ denote the {\em trim part} of $\A''$ (by removing all vertices/edges which are not part of successful
paths in $\A''$; this can be done effectively). Then $\A \sqcap L$ can be obtained by projecting
into the first component the various constituents of $\B$. \qed

Given an alphabet $A$, we denote by $A^{-1}$ a set of \emph{formal
  inverses} of $A$,
and write $\widetilde{A} = A\cup A^{-1}$. We say that $\widetilde{A}$
is an \emph{involutive alphabet}.
We extend $^{-1}:A \to A^{-1}: a \mapsto a^{-1}$ to an involution on
$\widetilde{A}^*$ through
\begin{displaymath}
(a^{-1})^{-1} = a,\quad (uv)^{-1} = v^{-1}u^{-1}\quad (a \in A,\;
u,v \in \widetilde{A}^*)\, .
\end{displaymath}

An automaton $\cal{A}$ over an involutive alphabet $\widetilde{A}$ is
\emph{involutive} if, whenever $(p,a,q)$
is an edge of $\cal{A}$, so
is $(q,a^{-1},p)$.  Therefore it suffices to depict just the
\emph{positively labelled} edges (having label in $A$) in their
graphical representation.

An involutive automaton is \emph{inverse} if it is deterministic, trim and has a single final state (note that for
involutive automata, being trim is equivalent to being connected). If
the latter happens to be the
initial state, it is called the \emph{basepoint}.




The next result is folklore. For a proof, see \cite[Proposition
2.2]{BS1}.

\bp
\label{invmor}
Given inverse automata $\cal{A}$ and $\cal{A}'$, then $L({\cal{A}})
\subseteq L({\cal{A}}')$ if and only if there exists a morphism
$\varphi: \cal{A} \to \cal{A}'$. Moreover, such a morphism is
unique.
\ep

Given an alphabet $A$, let $\sim$ denote the congruence on
$\widetilde{A}^*$ generated by the relation
\beq
\label{freegr}
\{(aa^{-1},1)\mid a \in \widetilde{A}\}\, .
\eeq
The quotient $F_A = \widetilde{A}^*/{\sim}$ is the \emph{free group
  on} $A$. We denote by $\theta: \widetilde{A}^* \to F_A$ the canonical
morphism $u \mapsto [u]_{\sim}$.

Alternatively, we can view (\ref{freegr}) as a {\em
  confluent}
length-reducing rewriting system on $\widetilde{A}^*$, where each word
$w \in \widetilde{A}^*$ can be transformed into a unique
\emph{reduced} word $\overline{w}$
with no factor of the form $aa^{-1}$.  As a
consequence, the equivalence
\begin{displaymath}
u\sim v \hspace{.5cm} \Leftrightarrow \hspace{.5cm}\overline{u} = \overline{v}
\hspace{2cm} (u,v \in \widetilde{A}^*)
\end{displaymath}
solves the word problem for $F_A$.
We shall use the notation
$R_A = \overline{\widetilde{A}^*}$.




We close this section with the classical Benois Theorem, which relates rational
languages with free group reduction:

\bt
\label{benois}
\cite{Ben}
If $L \subseteq \wt{A}^*$ is rational, then $\oo{L}$ is
an  effectively constructible rational language.
\et

\section{Stallings sections}

Let $G$ be a (finitely generated) group generated by the finite set $A$. More precisely, we consider an epimorphism
$\pi: \wt{A}^* \to G$ satisfying \beq \label{match} a\inv\pi = (a\pi)\inv \quad \mbox{ for every }\quad a \in A. \eeq A
homomorphism satisfying condition (\ref{match}) is said to be {\em
  matched}. Note that in this case (\ref{match}) holds for arbitrary
words.
For short, we shall refer to a matched epimorphism $\pi:
\wt{A}^* \to G$ (with
$A$ finite) as a {\em m-epi}.

We shall call a language $S \subseteq \wt{A}^*$ a {\em
section} (for $\pi$) if $S\pi = G$ and $S\inv = S$. For
every $X \subseteq G$, we write
$$S_X = X\pi\inv \cap S.$$

We say
that an effectively constructible rational section
$S \subseteq R_A$ is a {\em Stallings section} for $\pi$
if, for all $g,h \in G$:
\bi
\item[(S1)]
$S_g$ is an effectively constructible rational language;
\item[(S2)]
$S_{gh} \subseteq \oo{S_gS_h}$. \ei Note that (S2) yields immediately \beq \label{stwogen} S_{g_1\cdots g_n} \subseteq
\oo{S_{g_1}\cdots S_{g_n}} \eeq for all $g_1,\ldots, g_n \in G$. Moreover, in (S1) it suffices to consider $S_{a\pi}$
for $a \in A$. Indeed, by (\ref{stwogen}), and since $S\inv = S$ and $S_g\pi = g$ for every $g \in G$, we may write
$$S_{(a_1\cdots a_n)\pi} = \oo{S_{a_1\pi}\cdots S_{a_n\pi}} \cap S$$
and $S_{a_i\inv\pi} = S_{a_i\pi}\inv$ for all $a_i \in \wt{A}$. Then, by Proposition~\ref{proprat} and
Theorem~\ref{benois}, $S_g$ is a rational language for every $g\in G$; furthermore, it is effectively constructible
from $S_{a_1\pi},\ldots ,S_{a_n\pi}$.

Note that if $S$ is a Stallings section, then $S \cup \{ 1\}$ is
also a Stallings section. Indeed, it
is easy to see that conditions (S1) and (S2) are still
verified: namely, if $gh = 1$, then $1 \in \oo{S_gS_g\inv}
= \oo{S_gS_h}$ and so $S_{gh} \cup \{ 1 \} \subseteq \oo{S_gS_h}$ as
required.

The next result shows that the existence of a Stallings section is
independent from
the finite set $A$ and the m-epi $\pi:\wt{A}^* \to G$ considered:

\bp
\label{independence}
Let $\pi:\wt{A}^* \to G$ and $\pi':\wt{A'}^* \to G$ be m-epis. Then $G$
has a Stallings section for $\pi$ if and only if
$G$ has a Stallings section for $\pi'$.
\ep

\proof
Let $S \subseteq R_A$ be a Stallings section for $\pi$.
There exists a matched homomorphism $\p:\wt{A}^* \to \wt{A'}^*$
such that $\p\pi' = \pi$. Write $S' = \oo{S\p}$.
By Proposition \ref{proprat}(ii) and Theorem \ref{benois}, $S'$ is an
effectively constructible rational subset of $R_{A'}$.
We claim that
\beq
\label{independence1}
S'_g = \oo{S_g\p}
\eeq
holds for every $g \in G$.

Indeed, let $u \in S'_g$. Then $u = \oo{v\p}$ for some $v \in S$ and
$v\pi = v\p\pi' = \oo{v\p}\pi' = u\pi' = g$. Hence $v \in S_g$ and so
$S'_g \subseteq \oo{S_g\p}$.

Conversely, let $v \in S_g$. Then $\oo{v\p} \in \oo{S\p} = S'$ and
$\oo{v\p}\pi' = v\p\pi' = v\pi = g$, hence $\oo{v\p} \in S'_g$ and so
(\ref{independence1}) holds.

Since
$$(S')\inv = (\oo{S\p})\inv = \oo{(S\p)\inv} = \oo{S\inv\p} = \oo{S\p} = S',$$
it follows from
(\ref{independence1}) that $S'$ is a section for $\pi'$. Moreover,
(S1) is inherited by $S'$ from $S$ by Proposition
\ref{proprat}(ii) and Theorem \ref{benois}. Finally, for all
$g,h \in G$, we get
$$\begin{array}{lll}
S'_{gh}&=&\oo{S_{gh}\p} \subseteq \oo{(\oo{S_{g}S_h})\p} =
\oo{(S_{g}S_h)\p}\\
&=&\oo{(S_{g}\p)(S_h\p)} = \oo{(\oo{S_{g}\p})(\oo{S_h\p})} =
\oo{S'_{g}S'_h},
\end{array}$$
hence (S2) holds for $S'$ and so $S'$ is a Stallings section for $\pi'$.
By symmetry, we get the required equivalence.
\qed


\bp
\label{freefin}
Free groups of finite rank and finite groups have Stallings sections.
\ep

\proof
Let $A$ be a finite set and consider the canonical m-epi
$\theta: \wt{A}^* \to F_A$. Let $S = R_A = \oo{\wt{A}^*}$, which is
rational by Theorem \ref{benois}. Since $S_g = \oo{g}$ for
every $g \in F_A$, it is immediate that $S$ is a Stallings
section for $\theta$.

Assume now that $G$ is finite and $\pi: \wt{A}^* \to G$ is a m-epi. We
show that $S = R_A$ is a Stallings
section for $\pi$. For every $g \in G$, we have $S_g = g\pi\inv \cap
R_A = \oo{g\pi\inv}$. Since both $g\pi\inv$ and
$R_A$
are effectively
constructible
rational languages, so is their intersection and so (S1)
holds. Finally, let
$u \in S_{gh}$ and take $v \in S_h$. Then $(uv\inv)\pi = ghh\inv = g$
and so $\oo{uv\inv} \in \oo{g\pi\inv} = S_g$. Hence
$u = \oo{uv\inv v} = \oo{\oo{uv\inv} v} \in \oo{S_gS_h}$ and (S2)
holds as well. Therefore $R_A$ is a Stallings
section for $\pi$.
\qed

Given a m-epi $\pi:\wt{A}^* \to G$ and $H \leqslant G$, we define the
{\em Schreier automaton}
$\Gamma(G,H,\pi)$ to be the $\wt{A}$-automaton having:
\bi
\item
the right cosets $Hg$ $(g \in G)$ as vertices;
\item
$H$ as the basepoint;
\item
edges $Hg \mapright{a} Hg(a\pi)$ for all $g \in G$ and $a \in \wt{A}$.
\ei
It is immediate that $\Gamma(G,H,\pi)$ is always an inverse
$\wt{A}$-automaton, but it is infinite unless $H$ has finite index in
$G$. Moreover, $L(\Gamma(G,H,\pi)) = H\pi\inv$.

We will prove that $\Gamma(G,H,\pi) \sqcap S$ is an effectively
constructible finite inverse automaton when $S$ is a Stallings section for
$\pi$. The following lemmas pave the way for the construction of
$\Gamma(G,H,\pi) \sqcap S$:

\bl
\label{addinv}
Let $\pi: \wt{A}^* \to G$ be a m-epi. Let $\A$ be a trim
$\wt{A}$-automaton and let $p \mapright{a} q$ be an
edge of $\A$ for some $a \in \wt{A}$. Let $\B$ be obtained by adding
the edge $q \mapright{a\inv} p$ to $\A$. Then $(L(\B))\pi \subseteq
\langle (L(\A))\pi \rangle$.
\el

\proof Write $\A = (Q,q_0,T,E)$. We can factor any $u \in L(\B)$ as $u = u_0a\inv u_1 \cdots a\inv u_n$, where the
$a\inv$ label each visit to the new edge. We show that $u\pi \in \langle (L(\A))\pi \rangle$ by induction on $n$. The
case $n = 0$ being trivial, assume that $n \geq 1$ and the claim holds for $n-1$. Writing $v = u_0a\inv u_1 \cdots
a\inv u_{n-1}$, we have a path in $\B$ of the form
$$q_0 \mapright{v} q \mapright{a\inv} p \mapright{u_n} t \in T.$$
Since $\A$ is trim, we have also a path
$$q_0 \mapright{w} p \mapright{a} q \mapright{z} t' \in T$$
in $\A$. By the induction hypothesis, we get $(vz)\pi \in \langle
(L(\A))\pi \rangle$ and so
$$u\pi = (va\inv u_n)\pi = ((vz)(z\inv a\inv w\inv)(wu_n))\pi
\in \langle (L(\A))\pi \rangle$$
as claimed.
\qed

\bl
\label{basepoint}
Let $\pi: \wt{A}^* \to G$ be a m-epi.
Let $\A = (Q,q_0,T,E)$ be a trim $\wt{A}$-automaton and let $\B$ be obtained by
identifying $q_0$ with some $t \in T$. Then
$(L(\B))\pi \subseteq \langle (L(\A))\pi \rangle$.
\el

\proof Let $u \in L(\B)$. We can factor it as $u = u_1\cdots u_n$, where $p_i \mapright{u_i} q_i$ is a path in $\A$
with $p_i,q_i \in \{ q_0,t \}$ $(i = 1,\ldots,n)$. In any case, there exist paths
$$q_0
\mapright{v_i} p_i, \quad q_i \mapright{w_i} t \in T$$
in $\A$ with $v_i, w_i \in L(\A) \cup \{ 1\}$. Since $v_iu_iw_i \in L(\A)$,
we get $u_i\pi = (v_i\inv(v_iu_iw_i)w_i\inv)\pi \in \langle
(L(\A))\pi \rangle$ for every $i$
and so $u\pi \in \langle (L(\A))\pi \rangle$ as well.
\qed

\bl \label{ident} Let $\pi:\wt{A}^* \to G$ be a m-epi. Let $\A$ be an involutive $\wt{A}$-automaton and let $p
\mapright{w} q$ be a path in $\A$ with $w\pi = 1$. Let $\B$ be obtained by identifying the vertices $p$ and $q$. Then
$L(\A) \subseteq L(\B)$ and $(L(\B))\pi = (L(\A))\pi$. \el

\proof The first inclusion is clear. Since $\A$ is involutive, we have also a path $q \mapright{w\inv} p$ in $\A$ and
$w\inv\pi = 1$. Clearly, every $u \in L(\B)$ can be lifted to some $v \in L(\A)$ by inserting finitely many occurrences
of the words $w,w\inv$, that is, we can get factorizations
$$u = u_0u_1\cdots u_n \in L(\B),\quad v =
u_0w^{\varepsilon_1}u_1\cdots w^{\varepsilon_n} u_n \in L(\A)$$
with $\varepsilon_1, \ldots, \varepsilon_n \in \{ -1,\; 1 \}$.
Since $u\pi = v\pi$, it follows that $(L(\B))\pi \subseteq
(L(\A))\pi$. The opposite inclusion holds trivially.
\qed

Since $(aa\inv)\pi = 1$ for every $a \in \wt{A}$, this same argument
proves that:

\bl \label{fold} Let $\pi:\wt{A}^* \to G$ be a m-epi. Let $\A$ be a finite involutive $\wt{A}$-automaton and let $\B$
be obtained by successively folding pairs of edges in $\A$. Then $L(\A) \subseteq L(\B)$ and $(L(\B))\pi = (L(\A))\pi$.
\el

The next lemma reveals how the automaton $\Gamma(G,H,\pi) \sqcap S$
can be recognized.

\bl
\label{theeq}
Let $S\subseteq R_A$ be a Stallings
section for the m-epi
$\pi:\wt{A}^* \to G$ and let $H \leqslant_{f.g.} G$. Let
$\A$ be a finite inverse
$\wt{A}$-automaton with a basepoint such that
\beq
\label{theeq1}
S_H \subseteq L(\A) \subseteq H\pi\inv,
\eeq
\beq
\label{theeq2}
\mbox{there is no path $p \mapright{w} q$ in $\A$ with $p\neq q$ and
  $w\pi = 1$.}
\eeq
Then $\Gamma(G,H,\pi) \sqcap S \cong \A \sqcap S$.
\el

\proof
Since $\A$ and $\Gamma = \Gamma(G,H,\pi)$ are both inverse automata with a
basepoint, and $L(\A) \subseteq H\pi\inv = L(\Gamma)$, it
follows from Proposition \ref{invmor} that there exists a morphism
$\p:\A \to \Gamma$. Suppose that
$p\p = q\p$ for some vertices $p,q$ in $\A$.
Take geodesics
$$q_0 \mapright{u} p,\quad q_0 \mapright{v} q$$
in $\A$, where $q_0$ denotes the basepoint. Since $p\p = q\p$, we have
$uv\inv \in L(\Gamma) = H\pi\inv$. Let $s_0 \in S_{(uv\inv)\pi}
\subseteq S_H$. Then
$s_0 \in L(\A)$ by (\ref{theeq1}) and so there is a path $p
\vlongmapright{u\inv s_0v} q$ in $\A$. Since $(u\inv s_0v)\pi = (u\inv
uv\inv v)\pi = 1$, it follows from (\ref{theeq2}) that $p = q$. Thus
$\p$ is injective.

It is immediate that $\p$ restricts to an injective morphism $\p':\A
\sqcap S \to \Gamma \sqcap S$. It remains to show that every edge of
$\Gamma \sqcap S$ is induced by some edge of $\A \sqcap S$.
Assume that $H \mapright{s} H$ is a
(successful) path in $\Gamma$ with $s \in S$. By (\ref{theeq1}), we
have $s \in L(\A)$ and the path $q_0 \mapright{s} q_0$ is mapped by
$\p'$ onto $H \mapright{s} H$. Since every edge of $\Gamma \sqcap S$
occurs in some path $H \mapright{s} H$, it follows that $\p'$ is an
isomorphism.
\qed

\bl
\label{decident}
Let $S\subseteq R_A$ be a Stallings
section for the m-epi
$\pi:\wt{A}^* \to G$ and let $H \leqslant_{f.g.} G$. Let
$\A$ be a finite inverse
$\wt{A}$-automaton with a basepoint such that
$S_H \subseteq L(\A) \subseteq H\pi\inv$. It is decidable, given two distinct
vertices $p,q$ of $\A$, whether or not there is some path $p
\mapright{w} q$ in $\A$ with $w\pi = 1$.
\el

\proof
Let $p,q$ be distinct vertices of $\A$ and let $q_0$ denote its basepoint. Take
geodesics $q_0 \mapright{u} p$ and $q_0
\mapright{v} q$, and let $s \in S_{(uv\inv)\pi}$. We claim that
there is a path $p \mapright{w} q$ in $\A$ with
$w\pi = 1$ if and only if $s \in L(\A)$.

Indeed, assume that $p \mapright{w} q$ is such a path. Then $uwv\inv
\in L(\A)$ and so $s\pi = (uv\inv)\pi = (uwv\inv)\pi \in H$. Thus
$s \in S_H \subseteq L(\A)$.

Conversely, assume that $s \in L(\A)$. Then there is a path $p
\vlongmapright{u\inv sv} q$ in $\A$. Since $(u\inv sv)\pi = (u\inv
uv\inv v)\pi = 1$, the lemma is proved.
\qed

\bt
\label{ecsch}
Let $S\subseteq R_A$ be a Stallings
section for the m-epi $\pi:\wt{A}^* \to G$ and let $H \leqslant_{f.g.} G$. Then
$\Gamma(G,H,\pi) \sqcap S$ is an effectively constructible finite inverse
$\wt{A}$-automaton with a basepoint such that
\beq
\label{ecsch1}
S_H \subseteq L(\Gamma(G,H,\pi) \sqcap S) \subseteq H\pi\inv.
\eeq
\et

\proof
Assume that $H = \langle h_1,\ldots, h_m \rangle$. For $i =
1,\ldots,m$, let $\A_i = (Q_i, q_i,t_i,E_i)$ be a finite trim
$\wt{A}$-automaton with a single initial  and a single terminal
vertex satisfying
\beq
\label{ulti}
S_{h_i} \subseteq \oo{L(\A_i)} \subseteq h_i\pi\inv
\eeq
(in the next section we shall discuss how to define such an automaton
with the lowest possible complexity).
Let $\B_0$ be the
$\wt{A}$-automaton obtained by taking the disjoint union of the $\A_i$
and then identifying all the $q_i$ into a single initial vertex $q_0$.

Suppose that $q_i \mapright{u} q_i$ is a path in $\A_i$. Take $v \in L(\A_i)$. Then $uv \in L(\A_i) \subseteq
h_i\pi\inv$ and so $u\pi = (uvv\inv)\pi = h_ih_i\inv = 1$. It follows easily that $(L(\B_0))\pi \subseteq (S_{h_1} \cup
\cdots \cup S_{h_m})\pi \subseteq H$.

Let $\B_1$ be the finite trim involutive
$\wt{A}$-automaton obtained from $\B_0$ by
adjoining edges $(q,a^{-1},p)$ for all edges $(p,a,q)$ in $\B_0$ $(a \in
\wt{A})$. It follows from Lemma \ref{addinv} that $(L(\B_1))\pi \subseteq
\langle (L(\B_0))\pi \rangle \subseteq H$.

Next let $\B_2$ be the $\wt{A}$-automaton obtained from $\B_1$ by
identifying all terminal vertices with the initial vertex
$q_0$. By Lemma \ref{basepoint}, we get $(L(\B_2))\pi \subseteq
\langle (L(\B_1))\pi \rangle \subseteq H$.

Finally, let $\B_3$ be the finite inverse $\wt{A}$-automaton with a
basepoint obtained by
complete folding of $\B_2$. By Lemma \ref{fold}, we have $(L(\B_3))\pi
= (L(\B_2))\pi \subseteq H$ and so $L(\B_3) \subseteq H\pi\inv$.
Moreover,
$$S_{h_1} \cup \cdots \cup S_{h_m} \subseteq \oo{L(\A_1)}
\cup \cdots \cup \oo{L(\A_m)} \subseteq \oo{L(\B_0)} \subseteq \oo{L(\B_3)}$$
and $S\inv = S$ yield
$$\oo{(S_{h_1} \cup \cdots \cup S_{h_m} \cup S_{h_1\inv} \cup \cdots \cup
S_{h_m\inv})^*} \subseteq \oo{L(\B_3)}$$
since $\B_3$ is involutive and has a basepoint, and therefore
$$\oo{(S_{h_1} \cup \cdots \cup S_{h_m} \cup S_{h_1\inv} \cup \cdots \cup
S_{h_m\inv})^*} \subseteq L(\B_3)$$
since $\B_3$ is inverse (the language of an inverse automaton is
closed under reduction since a word $aa\inv$ must label only
loops). In view of (\ref{stwogen}), it follows that
$S_h \subseteq L(\B_3)$ for every $h \in H$ and so $S_H
\subseteq L(\B_3)$. Therefore (\ref{theeq1}) holds for
$\B_3$.

However, (\ref{theeq2}) may not hold. Assume that the vertex set $Q'$
of $\B_3$ is totally
ordered. By Lemma \ref{decident}, we can
decide if that happens, and find all concrete instances
$$J = \{ (p,q) \in Q' \times Q' \mid p < q \mbox{ and there is some path $p
\mapright{w} q$ in $\B_3$ with }w\pi = 1\}.$$
Let $\B_4$ be the finite inverse
$\wt{A}$-automaton with a basepoint obtained by identifying all pairs
of vertices in $J$ followed by complete folding. Since the existence
of a path with label in $1\pi\inv$ is preserved through the
identification process, it follows from Lemmas \ref{ident} and
\ref{fold} that $\B_4$ still satisfies (\ref{theeq1}).

Suppose that there exists a path $p' \mapright{w'} q'$ in $\B_4$ with
$p' \neq q'$ and $w'\pi = 1$. We can lift $p'$ and $q'$ to vertices $p$
and $q$ in $\B_3$, respectively. It is straightforward to check that
the path $p' \mapright{w'} q'$ can be lifted to a path $p \mapright{w}
q$ in $\B_3$ by successively inserting in $w'$ factors of the form:
\bi
\item
$aa\inv$ $(a \in \wt{A})$ (undoing the folding operations);
\item
$z \in 1\pi\inv$ (undoing the identification arising from $r \mapright{z} s$ ). \ei Since $w'\pi = w\pi$, it follows
that either $(p,q) \in J$ or $(q,p) \in J$, and so $p' = q'$, a contradiction. Therefore $\B_4$ satisfies
(\ref{theeq2}). Now the theorem follows from Proposition \ref{clsq} and Lemma \ref{theeq}. \qed

We call $\Gamma(G,H,\pi) \sqcap S$ the {\em Stallings automaton} of $H$ (for a
given Stallings section $S$). Note that $\Gamma(F_A,H,\theta) \sqcap
R_A$ is the classical Stallings automaton of $H \leq_{f.g.} F_A$ when
we take $R_A$ as Stallings section (for the canonical
m-epi $\theta$).

Stallings automata provide a natural decision procedure for the
generalized word problem:

\bc
\label{gwp}
Let $S\subseteq R_A$ be a Stallings
section for the m-epi $\pi:\wt{A}^* \to G$ and let $H \leqslant_{f.g.} G$.
Then the following conditions are equivalent for every $g \in
G$:
\bi
\item[(i)] $g \in H$;
\item[(ii)] $S_g \subseteq L(\Gamma(G,H,\pi) \sqcap S)$;
\item[(iii)] $S_g \cap L(\Gamma(G,H,\pi)
\sqcap S) \neq \emptyset$.
\ei
Therefore the generalized word problem is
decidable for $G$.
\ec

\proof
(i) $\Rw$ (ii). If $g \in H$, then $S_g
\subseteq S_H \subseteq L(\Gamma(G,H,\pi) \sqcap S)$ by Theorem
\ref{ecsch}.

(ii) $\Rw$ (iii). Immediate since $S_g \neq \emptyset$ due to $S$
being a section.

(iii) $\Rw$ (i). Since $S_g \cap L(\Gamma(G,H,\pi)
\sqcap S) \subseteq g\pi\inv \cap H\pi\inv$.

Now decidability follows from (S1) and Theorem \ref{ecsch}.
\qed

We can also prove the following generalization of Benois Theorem:

\bt
\label{newben}
Let $S\subseteq R_A$ be a Stallings
section for the m-epi $\pi:\wt{A}^* \to G$ and let $L \subseteq \wt{A}^*$ be
rational. Then $S_{L\pi}$ is an effectively
constructible rational language.
\et

\proof Let $\p:\wt{A}^* \to \rat \wt{A}^*$ be the rational substitution defined by $a\p = S_{a\pi}$, for $a \in \wt{A}$
(note that $1\p =\{ 1\}$ and, for $u=a_1\cdots a_n$ ($a_i\in \wt{A}$), $u\p$ is not $S_{u\pi}$ but just
$S_{a_1\pi}\cdots S_{a_n\pi}$). We claim that \beq \label{newben1} S_{u\pi} = S \cap \oo{u\p} \eeq holds for every $u
\in L\setminus \{ 1\}$. Let $u = a_1\cdots a_n \in L$ $(a_i \in \wt{A})$. Then by (\ref{stwogen}) we get
$$S_{u\pi} = S_{(a_1\pi)\cdots(a_n\pi)} \subseteq \oo{S_{a_1\pi}\cdots
  S_{a_n\pi}} = \oo{(a_1\p)\cdots(a_n\p)} = \oo{u\p}$$
and so $S_{u\pi} \subseteq S \cap \oo{u\p}$.

Since $a\p\pi = S_{a\pi}\pi = a\pi$ holds for every $a \in \wt{A}$,
the inclusion $S \cap \oo{u\p} \subseteq S_{u\pi}$ follows from
$\oo{u\p}\pi = u\p\pi = u\pi$. Therefore (\ref{newben1}) holds.

Now it becomes clear that
$$S_{L\pi} = S \cap (\cup_{u\in L} \overline{u\p} )= S \cap \overline{L\varphi}$$
if $1 \notin L$ and
$$S_{L\pi} = (S \cap \overline{L\varphi}) \cup S_1$$
if $1 \in L$.

Now $L\varphi$ is an effectively constructible rational language by (S1) and Proposition~\ref{proprat}(ii), and so is
$\overline{L\varphi}$ by Theorem \ref{benois}. Since $S$ and $S_1$ are rational, it follows from Proposition
\ref{proprat}(i) that $S_{L\pi}$ is rational and effectively constructible.
\qed

A natural question to ask at this stage is if we can identify a Stallings automaton for a given Stallings section $S$.
In the classical case of a free group $F_A$ with $S=R_A$ this is an elementary thing to do: in this case, an
$\wt{A}$-automaton $\A$ is of the form $\Gamma(F_A,H,\pi) \sqcap R_A=\S_H$ for some $H\leqslant_{f.g.} F_A$ if and
only if $\A$ is inverse, has a basepoint, and has no vertex of outdegree one except possibly the basepoint.

\bp
\label{identSA}
Let $S\subseteq R_A$ be a Stallings
section for a m-epi $\pi:\wt{A}^* \to G$. It is decidable, given a
finite $\wt{A}$-automaton $\A$, whether or not $\A \cong
\Gamma(G,H,\pi) \sqcap S$ for some $H \leqslant_{f.g.} G$.
\ep

\proof We may assume that $\A$ is inverse and has a basepoint. Write $\A = (Q,q_0,q_0,E)$. The equality $\A = \A \sqcap
S$ is an obvious necessary condition, decidable by Lemma \ref{clsq}. Thus we may assume that $\A = \A \sqcap S$ (in
particular, $\A$ is trim).

Since $S\subseteq R_A$ and $\A$ is trim, it follows that only the basepoint may have outdegree 1,
and so $\A \cong \S(K) \cong \Gamma(F_A,K,\theta) \sqcap R_A$ for some $K \leqslant_{f.g.} F_A$ \cite[Proposition
2.12]{BS1}: the standard algorithm \cite[Proposition 2.6]{BS1} actually computes a finite subset $X \subseteq R_A$
projecting onto a basis $X\theta$ of $K$. Let $K' =\langle X\pi \rangle \leqslant_{f.g.} G$. We claim that $\A \cong
\Gamma(G,H,\pi) \sqcap S$ for some $H \leqslant_{f.g.} G$ if and only if $\A \cong \Gamma(G,K',\pi) \sqcap S$, a
decidable condition in view of Theorem \ref{ecsch}.

The converse implication being trivial, assume that $\A =
\Gamma(G,H,\pi) \sqcap S$ for some $H \leqslant_{f.g.} G$. Since words
of $1\pi\inv$ can only label loops in $\Gamma(G,H,\pi)$, it follows from
Lemma \ref{theeq} that we only need to show that
\beq
\label{theeq12}
S_{K'} \subseteq L(\A) \subseteq K'\pi\inv.
\eeq

Since $\A \cong \Gamma(F_A,K,\theta) \sqcap R_A$, it follows from
Theorem \ref{ecsch} that
$$X \subseteq R_A \cap K\theta\inv \subseteq L(\A) \subseteq K\theta\inv.$$
Since $K\theta\inv \subseteq K'\pi\inv$, we get $L(\A) \subseteq K'\pi\inv$. Finally, $X \subseteq L(\A) \subseteq
H\pi\inv$ yields $X\pi \subseteq H$ and so $K' \leqslant H$. Hence
$$S_{K'} \subseteq S_H \subseteq L(\A)$$ by
(\ref{ecsch1}) and so (\ref{theeq12}) holds. Thus $\A \cong
\Gamma(G,K',\pi) \sqcap S$ and we are done.
\qed

\section{Complexity}

In this section we discuss, for a given Stallings section, an
efficient way (from the viewpoint of
complexity) of constructing the automata $\A_i$ in the proof of
Theorem \ref{ecsch} and compute an upper bound for the complexity of
the construction of the Stallings automata $\Gamma(G,H,\pi) \sqcap S$.

We say that an $\wt{A}$-automaton is {\em uniterminal} if it has a single terminal vertex. It is easy to see that there
exist rational languages which fail to be recognized by any uniterminal automaton (e.g. $R_A$, since regular languages
recognizable by uniterminal automata and containing the empty word must have a basepoint and so they are
submonoids). However, we can prove the following:

\bl
\label{sist}
Let $S\subseteq R_A$ be a Stallings
section for the m-epi
$\pi:\wt{A}^* \to G$ and let $g \in G$. Then there exists a finite
trim uniterminal $\wt{A}$-automaton $\C_g$ satisfying
$$S_{g} \subseteq \oo{L(\C_g)} \subseteq g\pi\inv.$$
\el

\proof
Let $\C = (Q,i,T,E)$ be the minimum automaton of $S_g$ (or any other
finite trim automaton with a single initial vertex recognizing $S_g$)
and let $\C_g$
be obtained by identifying all the terminal vertices of $\C$. Clearly,
$\C_g$ is a finite trim uniterminal automaton and $S_{g} = L(\C) \subseteq
L(\C_g)$ yields $S_g = \oo{S_g} \subseteq \oo{L(\C_g)}$. It remains to
be proved that $(L(\C_g))\pi = g$.

Let $u \in L(\C_g)$. Then there exists a factorization $u = u_0u_1\cdots u_k$ such that
$$i \mapright{u_0} t_0,\quad s_1 \mapright{u_1} t_1,\quad \ldots,
\quad s_k \mapright{u_k} t_k$$ are paths in $\C$ with $s_j,t_j \in T$. Take a path $i \mapright{v_j} s_j$ in $\C$, for
$j=1,\ldots ,k$. Then $v_j, v_ju_j \in L(\C)$ and so $v_j\pi = (v_ju_j)\pi = g$. Hence $u_j\pi = (v_j\inv v_ju_j)\pi =
g\inv g = 1$ and so $u\pi = (u_0u_1\cdots u_k)\pi = u_0\pi = g$ since $u_0 \in L(\C) = S_g$. Thus $(L(\C_g))\pi = g$
and so $\oo{L(\C_g)} \subseteq g\pi\inv$ as required. \qed

We introduce next a multiplication of (finite trim) uniterminal automata:
given (finite trim) uniterminal
$\wt{A}$-automata $\A = (Q,i,t,E)$ and $\A' = (Q',i',t',E')$, let $\A
\ast \A' = (Q'',i,t',E'')$ be the (finite trim) uniterminal
$\wt{A}$-automaton obtained by taking the disjoint union of
the underlying graphs of $\A$ and $\A'$ and identifying $t$ with $i'$.

\bl
\label{multipl}
Let $S\subseteq R_A$ be a Stallings
section for the m-epi
$\pi:\wt{A}^* \to G$ and let $g,g' \in G$. Let $\A$ and $\A'$ be finite
trim uniterminal $\wt{A}$-automata satisfying
$$S_{g} \subseteq \oo{L(\A)} \subseteq g\pi\inv,\quad S_{g'} \subseteq
\oo{L(\A')} \subseteq g'\pi\inv.$$
Then
$$S_{gg'} \subseteq \oo{L(\A\ast\A')} \subseteq (gg')\pi\inv.$$
\el

\proof
Since $L(\A)L(\A') \subseteq L(\A\ast\A')$, we get in view of (S2)
$$S_{gg'} \subseteq \oo{S_gS_{g'}} \subseteq \oo{L(\A)L(\A')}
\subseteq \oo{L(\A\ast\A')}.$$

Now let $u \in L(\A\ast\A')$. Then $u$ labels a path in $\A\ast\A'$ of
the form
$$i \mapright{u_0} p \mapright{u_1} p \mapright{u_2} \,\, \cdots \,\,
\mapright{u_{k-1}} p \mapright{u_k} t',$$ where we emphasize all the occurrences of the vertex $p$ obtained through the
identification of $t$ and $i'$. Now it is easy to see that there exist paths $i \mapright{u_0} t$ in $\A$ and $i'
\mapright{u_k} t'$ in $\A'$. Moreover, for each $j = 1,\ldots,k-1$, there exists either a path $t \mapright{u_j} t$ in
$\A$ or a path $i' \mapright{u_j} i'$ in $\A'$. Now, in view of $(L(\A))\pi = g$ and $(L(\A'))\pi = g'$, we can use the
same argument as in the proof of Lemma \ref{sist} to show that $u_j\pi = 1$ for $j = 1,\ldots,k-1$. Hence $u\pi  =
(u_0u_1\cdots u_k)\pi = (u_0u_k)\pi = gg'$ and so $\oo{L(\A\ast\A')} \subseteq (gg')\pi\inv$ as required. \qed

In view of the preceding two lemmas, we can now set an algorithm to construct the automata $\A_i$ in the proof of
Theorem \ref{ecsch}. All we need for a start are the minimum automata of $S_{a\pi}$ for each $a \in A$ (or any other
finite trim automaton with a single initial vertex recognizing $S_{a\pi}$; this can be effectively constructed by
(S1)). Following the argument in the proof of Lemma~\ref{sist}, we may identify all the terminal vertices to get finite
trim uniterminal $\wt{A}$-automata $\C_{a\pi}$ satisfying
$$S_{a\pi} \subseteq \oo{L(\C_{a\pi})} \subseteq a\pi\pi\inv.$$
Note that, since $S\inv = S$, we get finite trim uniterminal
$\wt{A}$-automata $\C_{a\inv\pi}$
satisfying
$$S_{a\inv\pi} \subseteq \oo{L(\C_{a\inv\pi})} \subseteq a\inv\pi\pi\inv$$
by exchanging the initial and the terminal vertices in $\C_{a\pi}$ and
replacing each edge $p \mapright{b} q$ by an edge $q \mapright{b\inv}
p$.

Now, given $h_i \in G$, we may represent it by some reduced word $a_1\cdots a_n$ $(a_i \in \wt{A})$, and may compute
$$\A_i = ((\cdots (\C_{a_1\pi}\ast \C_{a_2\pi})\ast \C_{a_3\pi})\ast
\cdots )\ast \C_{a_n\pi}.$$
By Lemma \ref{multipl}, $\A_i$ is a finite trim uniterminal
$\wt{A}$-automaton satisfying
$$S_{h_i} \subseteq \oo{L(\A_i)} \subseteq h_i\pi\inv.$$

What is the maximum size of $\A_i$ relatively to $|h_i|$? What is the
time complexity of the algorithm for its construction? Note that we
start with only finitely many ``atomic'' automata $\C_{a\pi}$ $(a \in
A)$. Hence the number of vertices (edges) in $\A_i$ is a bounded
multiple of $|h_i|$, therefore is $O(|h_i|)$, and the time complexity
of the construction (disjoint union followed by identification of two
vertices, $|h_i|-1$ times) is also clearly $O(|h_i|)$. This is why we
gave ourselves (and the reader) the trouble of constructing the $\A_i$
this way instead of just taking the minimum automaton of $S_{h_i}$,
whatever that may be!

But what is the time complexity of the full algorithm leading to the
Stallings automaton $\Gamma(G,H,\pi) \sqcap S$? It is also useful to
discuss the complexity of the important intermediate $\B_3$ in the
proof of Theorem \ref{ecsch} since $\B_3$ suffices for such
applications as the generalized word problem: indeed, since $\B_3$
satisfies (\ref{theeq1}), we may replace $\Gamma(G,H,\pi) \sqcap S$ by
$\B_3$ in Corollary \ref{gwp}.

Let $n = |h_1|+\cdots + |h_m|$. It follows easily from our previous discussion of the time complexity of the
construction of the $\A_i$ that $\B_0$ (and therefore $\B_1$ and $\B_2$) can be constructed in time $O(n)$. Since we
get to $\B_3$ through complete folding, the complexity of constructing $\B_3$ is that of the classical Stallings
construction in the free group. Touikan proved in \cite{Tou} that such complexity is $O(n\log^* n)$, where $\log^* n$
denotes the least integer $k$ such that the $k$th iterate of the log function of $n$ is at most 1 (for most practical
purposes, $O(n\log^* n)$ is similar to $O(n)$). Therefore $\B_3$ can be constructed in time $O(n\log^* n)$.

We shall now discuss the complexity of the construction of the
Stallings automata:

\bt \label{complex} Let $S\subseteq R_A$ be a Stallings section for the m-epi $\pi:\wt{A}^* \to G$ and let $H = \langle
h_1,\ldots, h_m\rangle \leqslant_{f.g.} G$. Then $\Gamma(G,H,\pi) \sqcap S$ can be constructed in time $O(n^3\log^*
n)$, where $n = |h_1|+\cdots + |h_m|$. \et

\proof
We go back to the proof of Theorem \ref{ecsch}, starting at $\B_3$.

The number of vertices of $\B_3$ is $O(n)$ and therefore we have
$O(n^2)$ candidate pairs to $J$. For each one of these pairs, we must
decide whether or not they belong to $J$. This involves bounding the
complexity of the algorithm described in the proof of Lemma
\ref{decident}.

Let $p,q$ be distinct vertices of $\B_3$ and let $q_0$ denote its
basepoint. Take
geodesics $q_0 \mapright{u} p$ and $q_0
\mapright{v} q$.
Clearly, $g = (uv\inv)\pi$ can be represented by a word of length
$O(n)$. It follows from the previous
discussion
on the complexity of the construction of $\A_i$ that we may construct
a finite trim uniterminal $\wt{A}$-automaton $\C_{g}$ satisfying
$$S_{g} \subseteq \oo{L(\C_{g})} \subseteq g\pi\inv$$
in time $O(n)$. Performing a complete folding on $\C_g$ (in time
$O(n\log^* n)$), we get a finite inverse $\wt{A}$-automaton $\D_{g}$ satisfying
$$S_{g} \subseteq L(\D_{g}) \subseteq g\pi\inv.$$
Since $S$ is a constant for our problem, we can compute an element
$s \in S \cap L(\D_{g}) = S_g$ in time $O(n)$ and check if $s
\in L(\B_3)$ in time $O(n)$. Therefore, by  the proof of Lemma
\ref{decident}, we can decide whether or not $(p,q) \in J$ in time
$O(n\log^* n)$. Since we had $O(n^2)$ candidates to consider, we may
compute $J$ in time $O(n^3\log^* n)$. It is very likely that this
upper bound can be improved.

Since $\B_4$ is obtained from $\B_3$ by identifying the pairs in $J$
followed by complete folding, and $\B_3$ has $O(n)$ vertices, then
$\B_4$ can be constructed in time $O(n^3\log^* n)$ in view of Touikan's
bound.

For the last step, we must discuss the time complexity of the
algorithm in the proof of Proposition \ref{clsq}. Note that $\B_4$ has
$O(n)$ vertices and therefore (since the alphabet is fixed) $O(n)$
edges. Since $S$ is a constant for our problem, we can build the
direct product of $\B_4$ by some deterministic automaton recognizing
$S$ in time
$O(n)$ and compute its trim part in time $O(n)$ (we have $O(n)$
vertices and $O(n)$ edges), and the final projection can also be
performed in linear time. Therefore $\Gamma(G,H,\pi) \sqcap S$ can be
constructed in time $O(n^3\log^* n)$, which means very close to
cubic complexity.
\qed

We should stress that the above discussion of time complexity was
performed for a
fixed Stallings section of a fixed group. But the computation of a
Stallings section for a (virtually free) group can be in itself a
costly procedure, particularly if it is supported by Bass-Serre theory
as in the present case. This will become more evident throughout the
next two sections.

\section{Amalgamation over finite groups}

Given groups $H$, $G_1$ and $G_2$, and isomorphisms $\p_j:H \to H_j
\leq G_j$ $(j = 1,2)$, the {\em free product with amalgamation}
(amalgam for short) of $G_1$ and $G_2$, relative to $\p_1$ and $\p_2$,
is defined as the quotient of the free product $G_1 \ast G_2$ by the
normal subgroup generated by the elements of the form
$(h\p_1)(h\inv\p_2)$ $(h \in H)$. It is usually denoted by $G_1 \ast_H
G_2$, whenever a specific reference to the homomorphisms $\p_j$ can be
omitted.


The groups $G_j$ embed canonically into $G_1 \ast_H G_2$, and we shall
actually view them as subgroups of their amalgam. In particular, we
view $H_1 = H_2$ as a subgroup of $G = G_1 \ast_H G_2$.


A factorization $g = w_1\cdots w_n$ is said to be a {\em reduced form} for $g \in G_1 \ast_H G_2$ if: \bi
\item[(i)] $w_1 \in G_1 \cup G_2$;
\item[(ii)] $w_1 \notin H_1 \cup H_2$ if $n > 1$;
\item[(iii)] $w_i \in G_j \setminus H_{j} \Rightarrow w_{i+1} \in
  G_{j+1} \setminus H_{j+1}$
\ei
hold for all $i \in \{ 1,\ldots, n-1\}$ and $j \in \{ 1,2\}$ modulo 2.

Every element of $G_1 \ast_H G_2$ can be represented by a reduced
form, but the representation is
not in general unique. However, this representation can be strictly
controlled (see e.g \cite[Chapter IV]{LS}):

\bp \label{normalforms} Let $u = u_1\cdots u_m$ and $v = v_1\cdots v_n$ be reduced forms of $G_1 \ast_H G_2$. Then $u =
v$ holds in $G_1 \ast_H G_2$ if and only if one of the following conditions holds: \bi
\item[(i)] $m = n = 1$ and $u_1 = v_1 \in G_1 \cup G_2$;
\item[(ii)] $m = n = 1$ and $u_1 = h\p_j$, $v_1 = h\p_{j+1}$ for some
  $h \in H$ and $j \in \{ 1,2\}$ modulo 2;
\item[(iii)] $m = n > 1$ and there exist $z_1, \ldots, z_{n-1} \in H$
  and $j \in \{ 1,2\}$ modulo 2 such that
$$\begin{array}{l}
u_1 = v_1(z_1\p_j) \mbox{ in } G_j,\\
u_2 = (z_1\inv\p_{j+1})v_2(z_2\p_{j+1}) \mbox{ in } G_{j+1},\\
\hspace{.5cm}\cdots \\
u_{n-1} =
(z_{n-2}\inv\p_{j+n-2})v_{n-1}(z_{n-1}\p_{j+n-2}) \mbox{ in }
G_{j+n-2},\\
u_{n} = (z_{n-1}\inv\p_{j+n-1})v_{n} \mbox{ in } G_{j+n-1}.
\end{array}$$
\ei
\ep

The main theorem of this section is

\bt
\label{cuaofg}
Let $G_1$ and $G_2$ be groups with Stallings sections and let $H$ be a
finite group. Then $G_1 \ast_H G_2$ has also a Stallings
section.

\et

\proof
Let $S$ (respectively $T$) be a
Stallings section for the m-epi $\pi_1:\wt{A_1}^* \to G_1$ (respectively
$\pi_2:\wt{A_2}^* \to G_2$). We assume that
$\wt{A_1}^* \cap \wt{A_2}^* = 1$ and write $A = A_1 \cup A_2$. Let
$H$ be a finite group and consider
isomorphisms $\p_j:H \to H_j
\leq G_j$ $(j = 1,2)$. We denote by $G = G_1 \ast_H G_2$ the amalgam
of $G_1$ and $G_2$ relative to $\p_1$ and $\p_2$. Let $\pi: \wt{A}^*
\to G$ be the m-epi induced by $\pi_1$ and $\pi_2$.

Let $B = \{ b_h \mid h \in H \}$ be a new alphabet and let $\psi:B^*
\to H$ be the homomorphism defined by $b_h\psi = h$ $(h \in H)$. Let
$\xi:B^* \to \rat \wt{A}^*$ be the rational substitution defined by
$$b_h\xi = S_{h\p_1} \cup T_{h\p_2} \subseteq h\p_1\pi\inv = h\p_2\pi\inv.$$
We define
$$L = 1\psi\inv\xi.$$
In the next lemma, we collect some important properties of $L$:

\bl
\label{elle}
\bi
\item[(i)] $L$ is an effectively constructible rational language;
\item[(ii)] $1 \in L$ and $L\pi = 1$;
\item[(iii)] $L^2 = L = L\inv$;
\item[(iv)] $(b_h\xi)L(b_{h\inv}\xi) \subseteq L$ for every $h \in H$.
\ei
\el

\proof
(i) Since $H$ is finite, $1\psi\inv$ and $L$ are rational and effectively
constructible by Proposition \ref{proprat}.

(ii) Indeed, if $(b_{h_1}\cdots b_{h_n})\psi = 1$, then $h_1\cdots h_n = 1$ and so $$(b_{h_1}\cdots b_{h_n})\xi\pi
\subseteq ((h_1\p_1\pi\inv)\cdots (h_n\p_1\pi\inv))\pi = (h_1\cdots h_n)\p_1 = 1.$$

(iii) The equality $L^2 = L$ follows from $(1\psi\inv)^2 = 1\psi\inv$. Now let $u \in L$. We may write $u \in
(b_{h_1}\cdots b_{h_n})\xi$ with $(b_{h_1}\cdots b_{h_n})\psi = 1$. It follows that $(b_{h_n\inv}\cdots
b_{h_1\inv})\psi = 1$. Since $S\inv = S$ and $T\inv = T$, we get
$$b_{h\inv}\xi = S_{h\inv\p_1} \cup T_{h\inv\p_2} = S_{h\p_1}\inv \cup
T_{h\p_2}\inv = (b_{h}\xi)\inv$$
for every $h \in H$ and so
$$\begin{array}{lll}
u\inv&\in&((b_{h_1}\xi) \cdots (b_{h_n}\xi))\inv = (b_{h_n}\xi)\inv
\cdots (b_{h_1}\xi)\inv\\
&=&(b_{h_n\inv}\xi)
\cdots (b_{h_1\inv}\xi) = (b_{h_n\inv}
\cdots b_{h_1\inv})\xi \subseteq 1\psi\inv \xi = L.
\end{array}$$
Thus $L\inv \subseteq L$ and so also $L = (L\inv)\inv \subseteq
L\inv$. Therefore $L = L\inv$.

(iv) Assume that $(b_{h_1}\cdots b_{h_n})\psi = 1$ and $u \in (b_{h_1}\cdots b_{h_n})\xi$. Then $h_1\cdots h_n = 1$ and
so $hh_1\cdots h_nh\inv = 1$. It follows that
$$(b_h\xi)u(b_{h\inv}\xi) \subseteq
(b_hb_{h_1}\cdots b_{h_n}b_{h\inv})\xi \subseteq L.$$
\qed

Let $$S' = S \setminus \bigcup_{h \in H} S_{h\p_1}, \quad
T' = T \setminus \bigcup_{h \in H} T_{h\p_2}.$$
Since $H$ is finite and $S,T$ are both Stallings sections, then
$S',T'$ are both effectively constructible rational languages.
We define
\beq
\label{ssec}
V = \oo{LSL} \cup \oo{LTL} \cup
\oo{(1 \cup LS')(LT'LS')^*(L \cup LT'L)}.
\eeq
Since $L,S,T,S',T'$ are all effectively constructible rational
languages, so is $V$, in view of Proposition
\ref{proprat} and Theorem \ref{benois}. Since $S$ and $T$
are sections for $\pi_1$ and $\pi_2$, and $1 \in L$, it follows from the
representation of amalgams in reduced form that $V$ is a section for
$\pi$. In particular, note that $(S')\inv = S'$, $(T')\inv
= T'$ and so $V\inv = V$.
It remains to be proved that $V$ satisfies axioms
(S1) and (S2).


Now let
$g = g_1\cdots g_n$ be a reduced form of $G$. We claim that \beq \label{thekey} V_g = \oo{LW^{(1)}_{g_1}\cdots
LW^{(n)}_{g_n}L}, \eeq where $W^{(i)} = S$ if $g_i \in G_1$ and $W^{(i)} = T$ if $g_i \in G_2$. In particular, $V_g =
\oo{LS_{g}L} = \oo{LT_{g}L}$ if $g \in G_1 \cap G_2 = H_1 = H_2$.

We prove two cases, the others are similar:

\medskip

\noindent
\underline{Case} $n = 1$ and $g_1 \in G_1$:

\medskip

\noindent
We must prove that $V_g = \oo{LS_{g_1}L}$. Indeed, it is immediate that
$\oo{LS_{g_1}L} \subseteq V$ and $\oo{LS_{g_1}L}\pi = (LS_{g_1}L)\pi =
S_{g_1}\pi = g_1 = g$, hence $\oo{LS_{g_1}L} \subseteq
V_g$. Conversely, if $u \in V_g$, it is clear from Proposition
\ref{normalforms} and $L\pi = 1$ that we must have $u \in \oo{LSL} \cup
\oo{LTL}$. If $u \in \oo{LS_xL}$ for some $x \in G_1$, then $g_1 = g =
u\pi = x$ and so $u
\in \oo{LS_{g_1}L}$. Hence we assume that $u \in \oo{LT_yL}$ for some $y
\in G_2$. It follows
that $g_1 = g = u\pi = y$ and so $g \in G_1 \cap G_2 = H_1 = H_2$. We can
then write $g_1 = h\p_1$ and $y = h\p_2$ for some $h \in H$. By Lemma
\ref{elle}, we get
$$u \in \oo{LT_{h\p_2}L} \subseteq
\oo{LT_{h\p_2}S_{h\inv\p_1}S_{h\p_1}L} \subseteq \oo{L^2S_{h\p_1}L}
= \oo{LS_{g_1}L}$$
and so (\ref{thekey}) holds in this case.

\medskip

\noindent
\underline{Case} $n = 2k$ and $g_1 \in G_1\setminus H_1$:

\medskip

\noindent
We must prove that
$$V_g = \oo{LS_{g_1}LT_{g_2}\cdots
  LS_{g_{2k-1}}LT_{g_{2k}}L}.$$
Indeed, it is immediate that $\oo{LS_{g_1}LT_{g_2}\cdots
  LS_{g_{2k-1}}LT_{g_{2k}}L} \subseteq V$ and
$$\begin{array}{lll}
\oo{LS_{g_1}LT_{g_2}\cdots
  LS_{g_{2k-1}}LT_{g_{2k}}L}\pi&=&(LS_{g_1}LT_{g_2}\cdots
  LS_{g_{2k-1}}LT_{g_{2k}}L)\pi\\
&=&(S_{g_1}T_{g_2} \cdots S_{g_{2k-1}}T_{g_{2k}})\pi = g_1\cdots g_{2k}
= g,
\end{array}$$ hence $\oo{LS_{g_1}LT_{g_2}\cdots
  LS_{g_{2k-1}}LT_{g_{2k}}L} \subseteq
V_g$. Conversely, if $u \in V_g$, it is clear from Proposition \ref{normalforms} and $L\pi = 1$ that we must have $u
\in \oo{LS_{g'_1}LT_{g'_2}\cdots
  LS_{g'_{2k-1}}LT_{g'_{2k}}L}$ where
$$\begin{array}{rll}
g'_1&=&g_1(h_1\p_1),\\
g'_2&=&(h_1\inv\p_2)g_2(h_2\p_2),\\
&&\cdots\\
g'_{2k-1}&=&(h_{2k-2}\inv\p_1)g_{2k-1}(h_{2k-1}\p_1),\\
g'_{2k}&=&(h_{2k-1}\inv\p_2)g_{2k}
\end{array}$$
for some $h_1, \ldots, h_{2k-1} \in H$.
Since $S$ and $T$ satisfy (S2) and by Lemma \ref{elle}(iv), we get
$$\begin{array}{lll}
u&\in&\oo{LS_{g'_1}LT_{g'_2}\cdots
  LS_{g'_{2k-1}}LT_{g'_{2k}}L}\\
&\subseteq&
\oo{LS_{g_1}S_{h_1\p_1}LT_{h_1\inv\p_2}T_{g_2}\cdots
  T_{h_{2k-2}\p_2}LS_{h_{2k-2}\inv\p_1}S_{g_{2k-1}}
  S_{h_{2k-1}\p_1}LT_{h_{2k-1}\inv\p_2}T_{g_{2k}}L}\\
&\subseteq &\oo{LS_{g_1}LT_{g_2}\cdots
  LS_{g_{2k-1}}LT_{g_{2k}}L}
\end{array}$$
and so $V_g = \oo{LS_{g_1}LT_{g_2}\cdots
  LS_{g_{2k-1}}LT_{g_{2k}}L}$ as claimed.

The other cases are absolutely similar, therefore (\ref{thekey})
holds. Since $S$ and $T$ satisfy (S1), and by Proposition
\ref{proprat}, Theorem \ref{benois} and Lemma \ref{elle}(i), $V_g$ is
an effectively constructible rational language for
every $g \in G$. Therefore (S1) holds for $V$.

As a consequence of (\ref{thekey}), we have $$V_{g} = \oo{V_{g_1\cdots
    g_i}V_{g_{i+1}\cdots g_n}}$$ whenever $g = g_1\cdots g_n$ is a
reduced factorization. In particular, $V_{gg'} = \oo{V_gV_{g'}}$ holds if $g = g_1\cdots g_n$ and $g' = g'_1\cdots
g'_m$ are reduced factorizations with $g_n \in G_1\setminus H_1$ and $g'_1 \in G_2\setminus H_2$, or vice-versa. We
shall refer to this case as the favourable case.

Given $g \in G$, let $||g||$ denote the number $n$ of components in a reduced form $g_1\cdots g_n$ of $g$. Let $g,g'
\in G$. We prove that \beq \label{vstwo} V_{gg'} \subseteq \oo{V_gV_{g'}} \eeq by induction on $k = ||g||+||g'||$.


If $||g|| = ||g'|| = 1$, we may assume that $g,g' \in G_1$
or $g,g' \in G_2$, otherwise we have the favourable case and we are
done. Without loss of generality, we may assume that $g,g' \in
G_1$. Then (\ref{thekey}) and (S2) for $S$ yield
$$V_{gg'} = \oo{LS_{gg'}L} \subseteq \oo{LS_{g}S_{g'}L} \subseteq
\oo{LS_{g}L^2S_{g'}L} = \oo{V_gV_{g'}}.$$
Therefore (\ref{vstwo}) holds for $k = 2$.

Assume now that 
$||g||+||g'|| > 2$ and (\ref{vstwo}) holds for smaller values of $||g||+||g'||$. Let $g = g_1\cdots g_n$ and $g' =
g'_1\cdots g'_m$ be reduced decompositions of $g$ and $g'$.

We do not have to consider the favourable case, hence we may assume that $g_n, g'_1 \in G_1$ or $g_n, g'_1 \in G_2$. By
symmetry, we may assume that $g_n, g'_1 \in G_1$ and $n > 1$. Write $x = g_1\cdots g_{n-1}$ and $y = g_ng'_1\cdots
g'_m$. Then $||x||+||y||  < ||g||+||g'||$ and so the induction hypothesis yields
$$V_{gg'} = V_{xy} \subseteq \oo{V_xV_{y}}.$$
Suppose first that $m = 1$. Then (S2) for $S$ yields
$$V_y = \oo{LS_{g_ng'_1}L} \subseteq
\oo{LS_{g_n}S_{g'_1}L} \subseteq
\oo{LS_{g_n}L^2S_{g'_1}L} = \oo{LS_{g_n}LV_{g'}}$$
and so in view of (\ref{thekey}) we get
$$V_{gg'} \subseteq \oo{V_xV_{y}} \subseteq \oo{V_xLS_{g_n}LV_{g'}} =
\oo{V_gV_{g'}}.$$ Now suppose that $m > 1$. Then $V_y \subseteq \oo{V_{g_ng'_1}V_{g'_2\cdots g'_m}}$ by the induction
hypothesis and so
$$V_{gg'} \subseteq \oo{V_xV_{y}} \subseteq
\oo{V_xV_{g_ng'_1}V_{g'_2\cdots g'_m}} \subseteq
\oo{V_xV_{g_n}V_{g'_1}V_{g'_2\cdots g'_m}} = \oo{V_gV_{g'}}$$
by the favourable case.
Thus (\ref{vstwo}) holds and so (S2) holds for $V$. Therefore $V$
is a Stallings section for $\pi$ and the theorem is proved.
\qed

\section{HNN extensions over finite groups}

Given a subgroup $H$ of a group $K$ and a monomorphism $\p:H \to K$,
the {\em HNN extension} $HNN(K,H,\p)$ is the group defined by the
relative presentation
$$\langle K,t\mid tht\inv = h\p\; (h \in H) \rangle,$$
that is, is the quotient of the free product $K \ast F_{\{t\}}$ by the
normal subgroup generated by the elements of the form
$tht\inv(h\inv\p)$ $(h \in H)$. For details, the reader is referred to
\cite[Chapter IV]{LS}.

We use the standard notation $H_1 = H$ and $H_{-1} = H\p$. A given factorization of $g$, $g =
w_0t^{\varepsilon_1}w_1\cdots t^{\varepsilon_n}w_n$, is said to be a {\em reduced form} for $g \in HNN(K,H,\p)$ if: \bi
\item[(i)] $w_i \in K$;
\item[(ii)] $\varepsilon_i \in \{ -1, \; 1\}$;
\item[(iii)] $\varepsilon_{i+1} = -\varepsilon_i \Rightarrow w_{i}
  \notin H_{\varepsilon_i}$
\ei
hold for every possible $i$. In particular, $1$ is a reduced form.

Every element of $HNN(K,H,\p)$ can be represented by a reduced
form, but the representation is
not in general unique. However, this representation becomes clear as a
consequence of the classical Britton's Lemma, which we choose to
present in the following form:

\bp \label{britton} Let $g = u_0t^{\varepsilon_1}u_1\cdots t^{\varepsilon_n}u_n$ be a reduced form of $HNN(K,H,\p)$.
The alternative reduced forms for $g$ in $HNN(K,H,\p)$ are obtained by replacing each occurrence of $t$ by some element
of $\cup_{h \in H} (h\p)th\inv$. \ep


In particular, $1$ is the unique reduced form for the identity and so
both $K$ and $F_{\{ t\}}$ embed canonically into $HNN(K,H,\p)$.

\bt
\label{hnn}
Let $K$ be a group with a Stallings section and let $\p:H \to K$ be a
monomorphism for some finite
subgroup $H$ of $K$. Then $HNN(K,H,\p)$ has also a Stallings
section.
\et

\proof
Let $S$ be
a Stallings section for the m-epi $\eta:\wt{A}^* \to K$. Write $G =
HNN(K,H,\p)$, $B = A \cup \{
b \}$ and let $\pi:\wt{B}^* \to G$ be the m-epi defined by $a\pi
= a\eta$ $(a \in \wt{A})$ and $b\pi = t$.

Let $C = \{ c_h \mid h \in H \}$ be a new alphabet and let $\psi:C^*
\to H$ be the homomorphism defined by $c_h\psi = h$ $(h \in H)$. Let
$\xi:C^* \to \rat\wt{B}^*$ be the rational substitution defined by
$c_h\xi = S_{h} \cup b\inv S_{h\p}b$. We define
$$L = 1\psi\inv\xi = \{ (c_{h_1}\cdots c_{h_n})\xi \mid h_1\cdots h_n
= 1\}.$$
The next lemma summarizes some important properties of $L$:

\bl
\label{newelle}
\bi
\item[(i)] $L$ is an effectively constructible rational language;
\item[(ii)] $1 \in L$ and $L\pi = 1$;
\item[(iii)] $L^2 = L = L\inv$;
\item[(iv)] $(c_h\xi)L(c_{h\inv}\xi) \subseteq L$ for every $h \in H$;
\item[(v)] $L \subseteq ((H\eta\inv)b\inv (H\p\eta\inv)b)^*(H\eta\inv)$.
\ei
\el

\proof
Since
$$c_h\xi\pi = (S_{h} \cup b\inv S_{h\p}b)\pi = S_h\pi = S_h\eta =
h$$ for every $h
\in H$, the proof of Lemma \ref{elle} can be used with straightforward
adaptations to prove (i)-(iv).

On the other hand, since $1 \in H\eta\inv$ we have
$$L = 1\psi\inv\xi \subseteq C^*\xi \subseteq ((H\eta\inv) \cup b\inv
(H\p\eta\inv)b)^* \subseteq ((H\eta\inv)b\inv
(H\p\eta\inv)b)^*(H\eta\inv)$$
and so (v) holds as well.
\qed

Now let
$$N = (S\wt{b})^*S \setminus \wt{B}^*(bS_Hb\inv \cup b\inv
S_{H\p}b)\wt{B}^*$$ denote the set of all words in $(S\wt{b})^*S$
representing reduced forms of $G$. Let $\alpha:\wt{B}^* \to \rat\wt{B}^*$
be the rational substitution defined by $a\alpha = a$ $(a \in \wt{A})$
and $b\alpha = bL$, $b\inv\alpha = L\inv b\inv = Lb\inv$. We claim that
$$V = \oo{N\alpha}$$
is a Stallings section for $\pi$.

By Theorem \ref{newben} and Lemma \ref{newelle}, the languages
$S$, $S_H$, $S_{H\p}$ and $L$ are all rational and effectively
constructible. By Proposition \ref{proprat} and Theorem \ref{benois},
so are $N$, $N\alpha$ and $V$.
Since $N\pi = G$ and $1 \in L$, it follows that $V\pi = G$. Note that
$S\inv = S$ yields $N\inv
= N$, and together with $L\inv = L$, this yields $V\inv = V$. Thus $V$
is a section for $\pi$.

Let $g = u_0t^{\varepsilon_1}u_1\cdots t^{\varepsilon_n}u_n$ be a reduced form of $G$. We claim that \beq \label{hnn1}
V_g = \oo{(S_{u_0}b^{\varepsilon_1}S_{u_1}\cdots
  b^{\varepsilon_n}S_{u_n})\alpha}.
\eeq
Since $L\pi = 1$, we have $\alpha\pi = \pi$ and so
$$\oo{(S_{u_0}b^{\varepsilon_1}S_{u_1}\cdots
  b^{\varepsilon_n}S_{u_n})\alpha}\pi = (S_{u_0}b^{\varepsilon_1}S_{u_1}\cdots
  b^{\varepsilon_n}S_{u_n})\pi = u_0t^{\varepsilon_1}u_1\cdots
t^{\varepsilon_n}u_n = g,$$ hence $$\oo{(S_{u_0}b^{\varepsilon_1}S_{u_1}\cdots
  b^{\varepsilon_n}S_{u_n})\alpha} \subseteq V \cap g\pi\inv = V_g.$$

Conversely, let $w \in V_g$. Then there exists a reduced form $v_0t^{\delta_1}v_1\cdots t^{\delta_m}v_m$ of $G$ such
that $w \in \oo{(S_{v_0}b^{\delta_1}S_{v_1}\cdots
  b^{\delta_m}S_{v_m)}\alpha}$. Then
$$g = w\pi = v_0t^{\delta_1}v_1\cdots t^{\delta_m}v_m$$ and it follows
from Proposition \ref{britton} that $m = n$, $\delta_i = \varepsilon_i$ for $i = 1,\ldots, n$, and
$v_0t^{\varepsilon_1}v_1\cdots t^{\varepsilon_n}v_n$ can be obtained from $u_0t^{\varepsilon_1}u_1\cdots
t^{\varepsilon_n}u_n$ by replacing each occurrence of $t$ by some element of $\cup_{h \in H} (h\p)th\inv$. Assume that
$t^{\varepsilon_i}$ is replaced by $((h_i\p)th_i\inv)^{\varepsilon_i}$. For $i = 1,\ldots, n$, write
$$x_i = \left\{
\begin{array}{ll}
h_{i}&\mbox{ if }\varepsilon_i = -1\\
h_{i}\p&\mbox{ if }\varepsilon_i = 1
\end{array}
\right.
\quad
y_i = \left\{
\begin{array}{ll}
h_{i}&\mbox{ if }\varepsilon_i = 1\\
h_{i}\p&\mbox{ if }\varepsilon_i = -1
\end{array}
\right.$$
and also $y_0 = x_{n+1} = 1$. Then
$$v_i = y_i\inv u_ix_{i+1}$$
for  $i = 0,\ldots, n$.

Moreover, we claim that
\beq
\label{hnn2}
\oo{(S_{x_i}b^{\varepsilon_i}S_{y_i\inv})\alpha}
\subseteq
\oo{b^{\varepsilon_i}\alpha}.
\eeq
Assume that $\varepsilon_i =
1$. Then $$\oo{(S_{x_i}b^{\varepsilon_i}S_{y_i\inv})\alpha} =
\oo{S_{x_i}bLS_{y_i\inv}} = \oo{S_{h_{i}\p}bLS_{h_i\inv}}
= \oo{bb\inv S_{h_{i}\p}bLS_{h_i\inv}} \subseteq
\oo{b(c_{h_i}\xi)L(c_{h_i\inv}\xi)} \subseteq \oo{bL} =
  \oo{b^{\varepsilon_i}\alpha}$$
by Lemma \ref{newelle}(iv).

Similarly, (\ref{hnn2}) holds for $\varepsilon_i =
-1$.
Hence
$$\begin{array}{lll}
w&\in&\oo{(S_{v_0}b^{\varepsilon_1}S_{v_1}\cdots
  b^{\varepsilon_n}S_{v_n)})\alpha}\\
&=&\oo{(S_{u_0x_1}b^{\varepsilon_1}S_{y_1\inv u_1x_2}\cdots
  S_{y_{n-1}\inv u_{n-1}x_n} b^{\varepsilon_n}S_{y_n\inv u_n})\alpha}\\
&\subseteq&\oo{(S_{u_0}S_{x_1}b^{\varepsilon_1}S_{y_1\inv}S_{u_1}S_{x_2}\cdots
  S_{y_{n-1}\inv}S_{u_{n-1}}S_{x_n}
  b^{\varepsilon_n}S_{y_n\inv}S_{u_n})\alpha}\\
&\subseteq&\oo{(S_{u_0}b^{\varepsilon_1}S_{u_1}\cdots
  b^{\varepsilon_n}S_{u_n})\alpha}
\end{array}$$
and so (\ref{hnn1}) holds.

Since $K$ has a Stallings section, has decidable generalized word
problem and so we can effectively compute a reduced form for any given
element of $G$. Therefore (S1) follows from (\ref{hnn1}).

If $g \in G$ has a reduced form $g = w_0t^{\varepsilon_1}w_1\cdots t^{\varepsilon_n}w_n$, we write $||g|| = n$. We show
that $V_{gg'} \subseteq \oo{V_gV_{g'}}$ for all $g,g' \in G$ by induction on $||g||+||g'||$.

Let $g = w_0t^{\varepsilon_1}w_1\cdots t^{\varepsilon_n}w_n$ and $g' = w'_0t^{\varepsilon'_1}w'_1\cdots
t^{\varepsilon'_m}w'_m$ be reduced forms. If \beq \label{hnn3} gg' = w_0t^{\varepsilon_1}w_1\cdots
t^{\varepsilon_n}w_nw'_0t^{\varepsilon'_1}w'_1\cdots t^{\varepsilon'_m}w'_m \eeq is a reduced form, then
$$\begin{array}{lll}
V_{gg'}&=&\oo{(S_{w_0}b^{\varepsilon_1}S_{w_1}\cdots
  b^{\varepsilon_n}S_{w_nw'_0}b^{\varepsilon'_1}S_{w'_1}\cdots
  b^{\varepsilon'_m}S_{w'_m})\alpha}\\
&\subseteq&\oo{(S_{w_0}b^{\varepsilon_1}S_{w_1}\cdots
  b^{\varepsilon_n}S_{w_n}S_{w'_0}b^{\varepsilon'_1}S_{w'_1}\cdots
  b^{\varepsilon'_m}S_{w'_m})\alpha}\\
&\subseteq&\oo{(S_{w_0}b^{\varepsilon_1}S_{w_1}\cdots
  b^{\varepsilon_n}S_{w_n})\alpha(S_{w'_0}b^{\varepsilon'_1}S_{w'_1}\cdots
  b^{\varepsilon'_m}S_{w'_m})\alpha} = \oo{V_gV_{g'}},
\end{array}$$
hence
we may assume that (\ref{hnn3}) is not a reduced form and $V_{g_1g_2}
\subseteq \oo{V_{g_1}V_{g_2}}$ whenever $||g_1||+||g_2|| <
||g||+||g'||$.
In particular, $n,m > 0$ and either
\beq
\label{hnn4}
\varepsilon_n = -\varepsilon'_1 =
1, \quad w_nw'_0 \in H
\eeq
or
$$\varepsilon_n = -\varepsilon'_1 =
-1, \quad w_nw'_0 \in H\p.$$
The second case being analogous, we assume that (\ref{hnn4})
holds. Let $h = w_nw'_0$. Then
$$gg' =  (w_0t^{\varepsilon_1}w_1\cdots
t^{\varepsilon_{n-1}}w_{n-1})((h\p)w'_1t^{\varepsilon'_2}w'_2\cdots t^{\varepsilon'_m}w'_m).$$ Since $g_1 =
w_0t^{\varepsilon_1}w_1\cdots t^{\varepsilon_{n-1}}w_{n-1}$ and $g_2 = (h\p)w'_1t^{\varepsilon'_2}w'_2\cdots
t^{\varepsilon'_m}w'_m$ are both reduced forms, the induction hypothesis yields
$$\begin{array}{lll}
V_{gg'}&=&V_{g_1g_2} \subseteq \oo{V_{g_1}V_{g_2}} =
\oo{(S_{w_0}b^{\varepsilon_1}S_{w_1}\cdots
  b^{\varepsilon_{n-1}}S_{w_{n-1}})\alpha
  (S_{(h\p)w'_1}b^{\varepsilon'_2}S_{w'_2}\cdots
  b^{\varepsilon'_{m}}S_{w'_m})\alpha}\\
&\subseteq&\oo{(S_{w_0}b^{\varepsilon_1}S_{w_1}\cdots
  b^{\varepsilon_{n-1}}S_{w_{n-1}}S_{h\p}S_{w'_1}b^{\varepsilon'_2}S_{w'_2}\cdots
  b^{\varepsilon'_{m}}S_{w'_m})\alpha}.
\end{array}$$
Now
$$\oo{S_{h\p}\alpha} = S_{h\p} \subseteq \oo{b(b\inv
  S_{h\p}bS_{h\inv})S_hb\inv} \subseteq \oo{bLS_hLb\inv} = \oo{(bS_hb\inv)\alpha}
\subseteq (bS_{w_n}S_{w'_0}b\inv)\alpha$$ yields
$$\begin{array}{lll}
V_{gg'}&\subseteq&\oo{(S_{w_0}b^{\varepsilon_1}S_{w_1}\cdots
  b^{\varepsilon_{n-1}}S_{w_{n-1}}S_{h\p}S_{w'_1}b^{\varepsilon'_2}
  S_{w'_2}\cdots  b^{\varepsilon'_{m}}S_{w'_m})\alpha}\\
&\subseteq&\oo{(S_{w_0}b^{\varepsilon_1}S_{w_1}\cdots
  b^{\varepsilon_{n-1}}S_{w_{n-1}}b^{\varepsilon_n}S_{w_n}
  S_{w'_0}b^{\varepsilon'_1} S_{w'_1}b^{\varepsilon'_2}
  S_{w'_2}\cdots  b^{\varepsilon'_{m}}S_{w'_m})\alpha} =
  \oo{V_gV_{g'}}
\end{array}$$
and so (S2) holds. Therefore $V$ is a Stallings section for $\pi$.
\qed

\section{Virtually free groups}

Recall that a {\em pushdown $A$-automaton} is a sextuple of the form
$\A = (Q,q_0,T,D,d_0,\delta)$, where $Q$ and $D$ are finite
sets, $q_0 \in Q$, $T \subseteq Q$, $d_0 \in D$ and
$\delta$ is a finite subset of $$Q \times (A \cup 1) \times D \times Q
\times D^*.$$

A {\em configuration} of $\A$ is an element of $Q \times
D^*$. The pair $(q_0,d_0)$ is the {\em initial configuration}.
If $(q,a,d,p,u) \in \delta$, we write
$$(q,vd) \etra{a} (p,vu)$$ for every $v \in D^*$. We call
this relation an {\em elementary transition}. If we have a sequence
$$(q_0,w_0) \etra{a_1} (q_1,w_1) \etra{a_2} \cdots \etra{a_n}
(q_n,w_n)$$ for some $n \geq 0$, we write
$$(q_0,w_0) \vvlongtra{a_1\cdots a_n} (q_n,w_n)$$ and we refer to it
as a {\em transition}.
The language {\em accepted} by $\A$ (by final states) is defined by
$$L(\A) = \{ w \in A^* \mid (q_0,d_0) \tra{w} (t,u) \mbox{ for some $t \in
T$ and }u \in D^*\}.$$
A language $L \subseteq A^*$ is
{\em context-free} if $L = L(\A)$ for some pushdown
$A$-automaton $\A$. For details on pushdown automata, the reader is
referred to \cite[Chapter 6]{HMU}.

Recall that a group is {\em virtually free} if it has a free subgroup
of finite index. Some recent papers involving virtually free groups
include \cite{GHHR, LS3, LS2}.

We can now prove the main theorem of the paper:

\bt
\label{vfree}
A finitely generated group has a Stallings section if and only if it
is virtually free.
\et

\proof
It is known that finitely generated virtually free groups are, up to
isomorphism, the fundamental groups of graphs of groups where the
graph, the vertex groups and the edge groups are all finite \cite[Theorem
7.3]{SW}. Moreover,
they can be obtained from finite groups by finitely many successive
applications of free products with amalgamation over finite groups and
HNN extensions over finite groups \cite[Chapter 1, Example
3.5 (vi)]{DD}. Since finite groups have Stallings
sections by Proposition \ref{freefin}, it follows from Theorems
\ref{cuaofg} and
\ref{hnn} that every finitely generated virtually free group has a
Stallings section.

Conversely, assume that $S$ is a Stallings section for the m-epi
$\pi:\wt{A}^* \to G$. We show that the {\em word problem
submonoid} $1\pi\inv$ is context-free. By Muller and Schupp's Theorem \cite{MS},
 this implies that $G$ is virtually free.

By the remark following the definition of Stallings section in
Section 3, we can assume that $1 \in S_1$.

For every $a \in \wt{A}$, let $\A^a = (Q^a,q_0^a,T^a,E^a)$ be a finite
automaton recognizing $S_{a\pi}$. We define a pushdown
$\wt{A}$-automaton $\A = (Q,q_0,t,D,d_0,\delta)$ by $Q = (\cup_{a \in
  \wt{A}} Q^a) \cup \{ q_0,t\}$, $D = \wt{A} \cup \{ d_0 \}$ and
$$\begin{array}{lll}\delta&=&\{ (q_0,1,d_0,t,1) \} \cup
\{ (q_0,a,d_0,q_0^a,d_0) \mid a \in
  \wt{A} \}\\
&\cup&\{ (p^a,1,d_0,q^a,d_0b) \mid (p^a,b,q^a) \in E^a, \; a,b \in \wt{A}
\}\\
&\cup&\{ (p^a,1,c,q^a,\oo{cb}) \mid (p^a,b,q^a) \in E^a, \; a,b,c \in \wt{A}
\}\\
&\cup&\{ (t^a,b,d,q_0^b,d) \mid t^a \in T^a, \; a,b \in \wt{A},\; d
\in D \}\\
&\cup&\{ (t^a,1,d_0,t,1) \mid  t^a \in T^a, \; a \in \wt{A} \}.
\end{array}$$
We shall prove that $1\pi\inv = L(\A)$. First of all we note that
$(p^a,b,q^a) \in E^a$ implies $$(p^{a},d_0v)
\etra{1} (q^{a},d_0\oo{vb})$$ for all $b \in \wt{A}$ and $v \in R_A$, hence
\beq
\label{vfree1}
\mbox{If $p^a \mapright{u} q^a$ is a path in $\A^a$, then $(p^{a},d_0v)
\tra{1} (q^{a},d_0\oo{vu})$}
\eeq
holds for all $a \in \wt{A}$ and $v \in R_A$.

Assume now that $a_1\cdots a_n \in 1\pi\inv$, with $a_1, \ldots, a_n \in \wt{A}$. We may assume that $n > 0$. Then $1
\in S_1 = S_{(a_1\cdots a_n)\pi} \subseteq \oo{S_{a_1\pi}\cdots S_{a_n\pi}}$ and so there exist $u_i \in S_{a_i\pi} =
L(\A^{a_i})$ such that $\oo{u_1\cdots u_n} = 1$. It follows from (\ref{vfree1}) that
$$(q_0^{a_i},d_0\oo{u_1\cdots u_{i-1}}) \tra{1}
(t^{a_i},d_0\oo{u_1\cdots u_{i}})$$
for some $t^{a_i} \in T^{a_i}$. Hence
$$\begin{array}{lll}
(q_0,d_0)&\etra{a_1}&(q_0^{a_1},d_0) \tra{1} (t^{a_1},d_0u_1)
\etra{a_2} (q_0^{a_2},d_0u_1) \tra{1} (t^{a_2},d_0\oo{u_1u_2})
\etra{a_3} \cdots\\ &&\\
&\etra{a_{n-1}}&(q_0^{a_{n-1}},d_0\oo{u_1\cdots u_{n-2}}) \tra{1}
(t^{a_{n-1}},d_0\oo{u_1\cdots u_{n-1}}) \etra{a_{n}}
(q_0^{a_{n}},d_0\oo{u_1\cdots u_{n-1}})\\ &&\\
&\tra{1}&(t^{a_{n}},d_0\oo{u_1\cdots u_{n}}) = (t^{a_{n}},d_0) \etra{1} (t,1)
\end{array}$$
and so $a_1\cdots a_n \in L(\A)$. Thus $1\pi\inv \subseteq L(\A)$.

Conversely, let $a_1\cdots a_n \in L(\A)$, with $a_1,\ldots, a_n \in \wt{A}$. We may assume that $n > 0$. It follows
easily that there exists a sequence of transitions of the form
$$\begin{array}{lll}
(q_0,d_0)&=&(p_0,d_0w_0) \etra{a_1} (q_0^{a_1},d_0w_0) \tra{1} (p_1,d_0w_1)
\etra{a_2} \cdots \etra{a_n} (q_0^{a_n},d_0w_{n-1}) \tra{1}
(p_n,d_0w_n)\\ &&\\
&=&(p_n,d_0) \etra{1} (t,1)
\end{array}$$
for some $w_0,\ldots, w_n \in R_A$. Now, for $i = 1,\ldots,n$, we must
have $p_i \in T^{a_i}$ and $w_i = \oo{w_{i-1}u_i}$ for some $u_i
\in L(\A^{a_i}) = S_{a_i\pi}$. Hence
$$1 = w_n = \oo{w_0u_1\cdots u_n} = \oo{u_1\cdots u_n} \in
\oo{S_{a_1\pi}\cdots S_{a_n\pi}}$$ and so $1 \in (S_{a_1\pi}\cdots S_{a_n\pi})\pi = (a_1\cdots a_n)\pi$. Thus $L(\A)
\subseteq 1\pi\inv$ and so $1\pi\inv = L(\A)$. Therefore $1\pi\inv$ is context-free and so $G$ is virtually free. \qed

\section{Sections with good properties}

Having established that finitely generated virtually free groups are
precisely the groups with a Stallings section, we have now the
possibility of imposing stronger conditions on their Stallings
sections, with the purpose of allowing further applications of the
Stallings automata $\Gamma(G,H,\pi) \sqcap S$.

The technique is simple. Suppose that:
\bi
\item
every finite group has a Stallings section with property $P$;
\item
if $G_1$ and $G_2$ have Stallings sections with property $P$ and $H$ is a
finite group, then $G_1 \ast_H G_2$ has also a Stallings
section with property $P$;
\item
if $K$ has a Stallings section with property $P$ and $H$ is a finite
subgroup of $K$, then $HNN(K,H,\p)$ has also a Stallings
section with property $P$.
\ei
Then, in view of \cite[Chapter 1, Example
3.5 (vi)]{DD}, every finitely generated virtually free group has a
Stallings section with property $P$.

A good example is given by the concept of extendable Stallings
section, which will turn out to be useful to
characterize finite index subgroups.

Let $S$ be a Stallings section for the m-epi $\pi:\wt{A}^* \to G$. We
say that $S$ is {\em extendable} if, for every $u \in S$, there exists
some $v \in R_A$ such that $uv^* \subseteq S$ and
\beq
\label{almost}
u \in \pref(S_{(uv^nu\inv)\pi}) \mbox{ for almost all } n \in \N.
\eeq

In order to prove the next result, we consider the following condition on a
Stallings section $S$ for $\pi:\wt{A}^* \to G$:
\bi
\item[(N)]
If $G$ is not torsion-free, then $S_1 \neq 1$.
\ei





\bp
\label{exten}
Every finitely generated virtually free group has an extendable
Stallings section.
\ep

\proof
In fact, we show that such a group always has an extendable
Stallings section satisfying condition (N).

Following the script previously described, we start by considering a
m-epi $\pi:\wt{A}^* \to G$ with $G$ finite. Let $S = R_A$ and take $v
= 1$ for every $u \in S$. Hence $uv^* \subseteq S$.
Since $S_g = \oo{g\pi\inv}$ for every $g \in
G$, we claim that $\pref(S_g) = R_A$:

Let $w \in R_A$ and take $a \in \wt{A}$
such that $wa \in R_A$. Since $G$ is finite, there exists some $m \in
\N$ such that every element of $G$ can be represented by some word of
length $< m$. In particular, there exists some $z \in R_A$ such that
$((a^{-m}w\inv)\pi)g = z\pi$ and $|z| < m$. Hence $(wa^mz)\pi = g$ and
so $\oo{wa^mz} \in \oo{g\pi\inv} = S_g$. Since $wa^m \in R_A$ and $|z|
< m$, we get $w \in \pref(S_g)$ and so $\pref(S_g) = R_A$.

Therefore (\ref{almost}) holds and so $R_A$ is an extendable Stallings
section for $\pi:\wt{A}^* \to G$ when $G$ is finite. Moreover, if $G$
is nontrivial, then $S_1 = \oo{1\pi\inv}$ contains nonempty words and
so condition (N) holds for $S$.

Next, assume that $S$ (respectively $T$) is an extendable
Stallings section for the m-epi $\pi_1:\wt{A_1}^* \to G_1$ (respectively
$\pi_2:\wt{A_2}^* \to G_2$), satisfying condition (N). We assume that
$\wt{A_1}^* \cap \wt{A_2}^* = 1$ and write $A = A_1 \cup A_2$. Let
$H$ be a finite group and consider
isomorphisms $\p_j:H \to H_j
\leq G_j$ $(j = 1,2)$. Let $G = G_1 \ast_H G_2$ be the amalgam
of $G_1$ and $G_2$ relative to $\p_1$ and $\p_2$, and let $\pi: \wt{A}^*
\to G$ be the m-epi induced by $\pi_1$ and $\pi_2$.
We may assume that $H_1 < G_1$ and $H_2 < G_2$, otherwise $G \cong
G_2$ or $G \cong G_1$.
We claim that the
Stallings section $V$ (for $\pi$) defined in the proof of Theorem
\ref{cuaofg} is also extendable. We use all the notation introduced in
that proof.

Let $u \in V$. Without loss of generality, we may assume that either $u = 1$ or the last letter of $u$ is in
$\wt{A_1}$. Let $u\pi = g_1\cdots g_m$ be a reduced form of $u\pi$.

Suppose first that $g_m \in G_1\setminus H_1$. Take $w \in T'$ and $z
\in S'$, and write $v = wz$. Then $uv^* \subseteq V$. We claim that
\beq
\label{aclaim}
u \in \pref(V_{(uv^nu\inv)\pi})\mbox{ if }n \geq \frac{m}{2} +1.
\eeq

Indeed, a simple induction on $\ell$ shows that if $x = x_1\cdots x_{k+\ell+1}$ and $y = y_1\cdots y_{\ell}$ are
reduced forms in $G$, then $xy$ has a reduced form $x_1\cdots x_kz_1\cdots z_r$:

We may assume that $y_1 \in G_1$. If $y_1 \in H_1$, then $\ell = 1$ and $xy = x_1\cdots x_{k+1}(x_{k+2}y_1)$ is a
reduced form and we are done. Hence we may assume that $y_1 , x_{k+\ell+1}\in G_1\setminus H_1$ and $x_{k+\ell+1}y_1 =
h\p_1$ for some $h \in H$. Then $(h\p_2\,y_2)y_3 \cdots y_{\ell}$ is a reduced form. If $\ell = 1$, then $xy =
x_1\cdots x_{k}(x_{k+1}(h\p_2))$ is a reduced form and we are done. If $\ell > 1$, we reduce the problem to the case
$\ell -1$ by considering the product $xy = (x_1\cdots x_{k+\ell})((h\p_2\, y_2)y_3 \cdots y_{\ell})$. Thus $xy$ has a
reduced form $x_1\cdots x_kz_1\cdots z_r$ as claimed.

In particular, taking $x = (uv^n)\pi$and $y = u\inv\pi$, it follows that $(uv^nu\inv)\pi$ has a reduced form $g_1\cdots
g_mw\cdots$ if $n \geq \frac{m}{2}+1$. Since $g_m \in G_1\setminus H_1$ and $uw$ is reduced, it follows easily from
(\ref{thekey}) that $V_{(uv^nu\inv)\pi}$ must contain some word $uw\cdots$ Thus (\ref{aclaim}) holds.

Suppose next that $g_m \in G_2\setminus H_2$. We can of course assume
that $u \neq 1$. Since  the last letter of $u$ is in $\wt{A_1}$, it follows from
(\ref{thekey}) that  $H \neq 1$. Since $S$ and $T$ satisfy condition
(N), it follows that there exist $z_1 \in S_1 \setminus \{ 1 \}$ and
$w_1 \in T_1 \setminus \{ 1 \}$. Take $w \in T'$ and $z
\in S'$, and write $v = w_1zwz_1$. Since $w_1,z_1 \in L$, we have
$uv^* \subseteq V$. Similarly to the preceding case, (\ref{aclaim})
holds.

Finally, we are left with the case $m = 1$ and $g_1 \in H_1$ (equal to
$H_2$ in $G$). Let $v = 1$. Then $uv^* \subseteq V$ trivially and
$(uv^nu\inv)\pi = 1$, hence it suffices to show that $u \in
\pref(V_1)$. By (\ref{thekey}) and Lemma \ref{elle}, we have $V_1 =
\oo{LS_1L} = \oo{L}$ and $u \in \oo{LS_{h\p_1}L}$ for some $h \in
H$. Let $w \in T_{h\inv\p_2}$. Then
$u \in \pref(uw)$ and $uw = \oo{uw} \in \oo{LS_{h\p_1}LT_{h\inv\p_2}}
\subseteq \oo{L} = V_1$ by Lemma \ref{elle}. Therefore $V$ is
extendable.

Suppose now that $V_1 = \{ 1 \}$. Since $V_1 = \oo{L}$, it follows
that $H$ is trivial and $S_1 = T_1 = \{ 1 \}$. Since $S$ and $T$
satisfy condition (N), it follows that $G_1$ and $G_2$ are
torsion-free and so $G$ is a free product of torsion-free groups,
hence torsion-free. Therefore $V$ satisfies condition (N).

Finally, assume that $S$ is
an extendable Stallings section for the m-epi $\eta:\wt{A}^* \to
K$ satisfying condition (N). Let $\p:H \to K$ be a
monomorphism for some finite
subgroup $H$ of $K$.
Write $G =
HNN(K,H,\p)$, $B = A \cup \{
b \}$ and let $\pi:\wt{B}^* \to G$ be the m-epi defined by $a\pi
= a\eta$ $(a \in \wt{A})$ and $b\pi = t$.

We claim that the
Stallings section $V$ (for $\pi$) defined in the proof of Theorem
\ref{hnn} is also extendable. We use all the notation introduced in
that proof.

We start by proving the following lemma:

\bl \label{parti} Let $k_0t^{\varepsilon_1}k_1\cdots t^{\varepsilon_m}k_m$ be a reduced form of $G$ with $m \geq 1$ and
let
$$\begin{array}{lll}
P&=&\{ z_0w_1b^{\varepsilon_1}z_1w_2b^{\varepsilon_2}z_2\cdots
w_mb^{\varepsilon_m}z_m \mid z_0 \in (k_0H_{-\varepsilon_{1}})\eta\inv,\\
&&z_i \in (H_{\varepsilon_i}k_iH_{-\varepsilon_{i+1}})\eta\inv\mbox{
  for }i = 1,\ldots,n-1, \; z_m \in (H_{\varepsilon_m}k_m)\eta\inv,\\
&&w_j \in (b^{\varepsilon_j}(H_{\varepsilon_j}\eta\inv)
b^{-\varepsilon_j}(H_{-\varepsilon_j}\eta\inv))^* \}.
\end{array}$$
Then $P$ is closed under (partial) free group reduction.
\el

\proof Let $u = z_0w_1b^{\varepsilon_1}z_1w_2b^{\varepsilon_2}z_2\cdots w_mb^{\varepsilon_m}z_m$ be an element of $P$
of the described form. Suppose first that $aa\inv$ is a factor of $u$ for some $a \in \wt{A}$. Then $aa\inv$ is either
a factor of some $z_i$ or a factor of some $w_j$, and it follows from the definitions that we may cancel $aa\inv$ and
remain inside $P$.

Thus we are left to discuss the case of cancellations involving the letter $b$. Suppose first that we have a factor
$b^{\varepsilon_j}b^{-\varepsilon_j}$ to cancel in $w_j$. Write $w_j = x_1\cdots x_r$, with $x_i =
b^{\varepsilon_j}x'_ib^{-\varepsilon_j}x''_i$ $(x'_i \in H_{\varepsilon_j}\eta\inv, x''_i \in
H_{-\varepsilon_j}\eta\inv)$, and assume that $x'_{\ell} = 1$. Cancelling our factor yields $x_1\cdots
x_{\ell-1}x''_{\ell}x_{\ell+1}\cdots x_r$. If $\ell > 1$, we can incorporate $x''_{\ell}$ into $x_{\ell-1}$ in view of
$(H_{-\varepsilon_j}\eta\inv)^2 \subseteq H_{-\varepsilon_j}\eta\inv$. If $\ell = 1$, we can incorporate $x''_{\ell}$
into $z_{j-1}$ by the same reason. The case of cancellations $b^{-\varepsilon_j}b^{\varepsilon_j}$ inside $w_j$ is
discussed similarly.

Suppose now that $b^{-\varepsilon_j}$ is the last letter of $w_j$ and cancels with its right neighbour
$b^{\varepsilon_j}$. Then we replace $w_jb^{\varepsilon_j} = x_1\cdots x_rb^{\varepsilon_j}$ by $x_1\cdots
x_{r-1}b^{\varepsilon_j}x'_r$. Since $x'_r$ can be absorbed by $z_j$, the claim holds also in this case.

Finally, we note that we can never have $z_j = 1$ when
$\varepsilon_j = -\varepsilon_{j+1}$: otherwise, we would get
$$1 = z_j \in (H_{\varepsilon_j}k_jH_{\varepsilon_{j}})\eta\inv$$ and
so $k_j \in H_{\varepsilon_j}$, impossible since $\varepsilon_j = -\varepsilon_{j+1}$ and
$k_0t^{\varepsilon_1}k_1\cdots t^{\varepsilon_m}k_m$ is a reduced form. \qed

Back to the proof of Proposition \ref{exten}, let $u \in V$. Assume that $u\pi = k_0t^{\varepsilon_1}k_1\cdots
t^{\varepsilon_m}k_m$ is a reduced form of $u\pi$.
If $m = 0$, then $u \in V_{u\pi} = S_{u\eta}$ and it follows that $uv^* \subseteq V$ for $v = b$. What if $m > 0$? Then
it is clear that $k_0t^{\varepsilon_1}k_1\cdots t^{\varepsilon_m}k_mt^{\varepsilon_mn}$ is a reduced form for every $n
\geq 0$. We claim that $ub^{\varepsilon_mn} \subseteq V$.

Indeed, let $P$ be defined as in
Lemma \ref{parti}. It follows from (\ref{hnn1}) that
$$u \in \oo{(S_{k_0}b^{\varepsilon_1}S_{k_1}\cdots
t^{\varepsilon_m}S_{k_m})\alpha}$$ and it is immediate that $S_{k_0}b^{\varepsilon_1}S_{k_1}\cdots
t^{\varepsilon_m}S_{k_m} \subseteq P$. Since
$$b\inv\alpha = Lb\inv \subseteq
(H\eta\inv)(b\inv (H\p\eta\inv)b(H\eta\inv))^*b\inv$$ by Lemma \ref{newelle}(v), and the factors $z_{j-1}$ may absorbe
factors from $H\eta\inv$ on the right when $\varepsilon_j = -1$, it follows that we may replace $b\inv$ by
$b\inv\alpha$ in $S_{k_0}b^{\varepsilon_1}S_{k_1}\cdots t^{\varepsilon_m}S_{k_m}$ and remain inside $P$.

Similarly,
$$b\alpha = bL \subseteq
b((H\eta\inv)b\inv (H\p\eta\inv)b)^*(H\eta\inv) = (b(H\eta\inv)b\inv (H\p\eta\inv))^*b(H\eta\inv)$$ and the factors
$z_{j}$ may absorbe factors from $H\eta\inv$ on the left when $\varepsilon_j = 1$, it follows that
$(S_{k_0}b^{\varepsilon_1}S_{k_1}\cdots b^{\varepsilon_m}S_{k_m})\alpha \subseteq P$. Hence
$\oo{(S_{k_0}b^{\varepsilon_1}S_{k_1}\cdots t^{\varepsilon_m}S_{k_m})\alpha} \subseteq P$ by Lemma \ref{parti} and so
$u \in P$.

As a consequence, we may write $u = xb^{\varepsilon_m}y$ with $y \in R_A$. Since $k_0t^{\varepsilon_1}k_1\cdots
t^{\varepsilon_m}k_mt^{\varepsilon_mn}$ is a reduced form for every $n \in \N$, it follows that
$$\oo{(S_{k_0}b^{\varepsilon_1}S_{k_1}\cdots
b^{\varepsilon_m}S_{k_m})\alpha b^{\varepsilon_mn}} \subseteq
\oo{(S_{k_0}b^{\varepsilon_1}S_{k_1}\cdots
b^{\varepsilon_m}S_{k_m}b^{\varepsilon_mn})\alpha} \subseteq V$$
for every $n \in \N$ and so $\oo{ub^{\varepsilon_mn}} \in V$. Since $u
= xb^{\varepsilon_m}y$, taking $v = b^{\varepsilon_m}$ we get $uv^n =
ub^{\varepsilon_mn} =
\oo{ub^{\varepsilon_mn}} \in V$ for every $n \in \N$ as claimed.

We continue now by showing that in any case
\beq
\label{hclaim}
u \in \pref(V_{(uv^nu\inv)\pi})\mbox{ if }n > m.
\eeq

Indeed, since $uv \in R_B$, it suffices to show that $(uv^nu\inv)\pi$ has a reduced form \linebreak
$k_0t^{\varepsilon_1}k_1\cdots t^{\varepsilon_m}k_mt^{\varepsilon_m}\cdots$ if $n > m$ ($k_0t\cdots$ if $m = 0$), and
we may use induction on $m$. The case $m = 0$ being obvious, assume that $m > 0$ and the claim holds for $m-1$. We
assume that $v = t$, the other case being analogous. If $k_0t^{\varepsilon_1}k_1\cdots t^{\varepsilon_m}k_mt^nk_m\inv
t^{-\varepsilon_m}k_{m-1}\inv\cdots t^{-\varepsilon_1}k_0\inv$ is not itself a reduced form, then $k_m\inv \in H$ and
$\varepsilon_m = 1$, hence we may write
$$(uv^nu\inv)\pi = k_0t^{\varepsilon_1}k_1\cdots
t^{\varepsilon_m}k_mt^{n-1}k'_mk_{m-1}\inv t^{-\varepsilon_{m-1}}k_{m-2}\inv\cdots t^{-\varepsilon_1}k_0\inv$$ for some
$k'_m \in H\p$. Since $n-1 >m-1$, the induction hypothesis applied to the product $k_0t^{\varepsilon_1}k_1\cdots
t^{\varepsilon_{m-1}}k_{m-1}(k'_m)\inv$ yields now the required result. Thus (\ref{hclaim}) holds and $V$ is
extendable.

It remains to show that $V$ satisfies condition (N).
Suppose that $V_1 = \{ 1 \}$. Since $V_1 = S_1$ and $S$ satisfies
condition (N), it follows that $K$ is torsion-free, and so
$H$ is trivial. Hence $G$ is the free product of $K$ by the infinite
cyclic group $F_{\{ t\}}$. Being a free product of torsion-free groups,
$G$ is itself torsion-free, therefore $V$ satisfies condition (N).
\qed

We can now derive the following application of the concept of
extendable Stallings section:

\bt
\label{findex}
Let $S$ be an extendable Stallings section for the m-epi $\pi: \wt{A}^*
\to G$ and let $H$ be a finitely generated subgroup of $G$. Then
the following conditions are equivalent:
\bi
\item[(i)]
$H$ has finite index in $G$;
\item[(ii)] $S \subseteq \pref(S_H)$;
\item[(iii)] every word of $S$ labels a path out of the basepoint of
  $\Gamma(G,H,\pi) \sqcap S$.
\ei
\et

\proof
(i) $\Rw$ (ii). Suppose that $u \in S \setminus \pref(S_H)$. Since $S$
is extendable, there exist some $v \in R_A$ and $m \in \N$ such that
$uv^* \subseteq S$ and $u \in \pref(S_{(uv^nu\inv)\pi})$ for $n \geq
m$. We claim that
\beq
\label{findex1}
H(uv^j)\pi \neq H(uv^i)\pi \mbox{ if } j \geq i+m.
\eeq
Indeed, assume that $j \geq i+m$. If $H(uv^j)\pi = H(uv^i)\pi$,
then $(uv^{j-i}u\inv)\pi \in H$ and so
$$u \in \pref(S_{(uv^{j-i}u\inv)\pi}) \subseteq \pref(S_H),$$
a contradiction. Therefore (\ref{findex1}) holds and so
$H$ has infinite index in $G$.

(ii) $\Rw$ (iii). Since $S_H \subseteq L(\Gamma(G,H,\pi) \sqcap S)$ by
Theorem \ref{ecsch}.

(iii) $\Rw$ (i). Assume now that every word of $S$ labels a path out
of the basepoint $q_0$ of
  $\A = \Gamma(G,H,\pi) \sqcap S$. Let $Q$ denote the vertex set of $\A$.
For every $q \in Q$, fix a path $q_0 \mapright{w_q}
q$. We claim that
\beq
\label{fid}
G = \bigcup_{q \in Q} H(w_q\pi).
\eeq
Indeed, let $g \in G$, and take $u \in S_g$. Then there is a path in
$\A$ of the form $q_0 \mapright{u}
q$ for some $q \in Q$. Hence $uw_q\inv \in L(\A) \subseteq H\pi\inv$
by Theorem \ref{ecsch} and so $g = u\pi \in H(w_q\pi)$. Thus
(\ref{fid}) holds and so $H$ has finite index in $G$.
\qed

A natural question to ask is whether or not one could replace
condition (S2) in the
definition of Stallings section by the stronger condition
\bi
\item[(S2')]
$S_{gh} = \oo{S_gS_h}$ for all $g,h \in G$.
\ei
However, we can prove that this condition can only be assumed in the
simplest cases:

\bp
\label{strong}
The following conditions are equivalent for a group $G$:
\bi
\item[(i)] there exist a m-epi $\pi:\wt{A}^* \to G$ and a Stallings
  section $S$ for $\pi$ satisfying (S2');
\item[(ii)] $G$ is either finite or free of finite rank;
\item[(iii)] $R_A$ is a Stallings
  section for some m-epi $\pi:\wt{A}^* \to G$.
\ei
\ep

\proof
(i) $\Rw$ (ii). Let $S$ be a Stallings
  section $S$ for $\pi:\wt{A}^* \to G$ satisfying (S2'). Then $S_1\inv
  = S_1 = \oo{S_1^2}$ and so we can view $(S_1,\circ)$ as a subgroup of
  $(R_A,\circ) \cong F_A$, where $u \circ v = \oo{uv}$. The same holds
  for $(S,\circ)$ since $S\inv
  = S = \oo{S^2}$, and $(S_1,\circ)$ is then a subgroup of
  $(S,\circ)$. Now $(S,\circ)$ must be free by Nielsen's
  Theorem. Since $S$, being a
  Stallings section, is rational, so is $(S,\circ)$ (a rational
  expression for $S$ as a subset of $\wt{A}^*$ translates through
  reduction to a rational
  expression for $S$ as a subset of $(R_A,\circ)$). The same happens
  with $S_1$, so it follows from Anisimov and Seifert's Theorem
  \cite[Theorem 3.1]{BS1}
  that both $(S,\circ)$ and $(S_1,\circ)$ are finitely generated
  groups.
Hence $(S,\circ)$ is a free group of finite rank.

For every $u \in S$, we have
$$\oo{uS_1u\inv} \subseteq \oo{S_{u\pi}S_1S_{u\inv\pi}} = S_1,$$
hence $(S_1,\circ)$ is a finitely generated normal subgroup of
the free group $(S,\circ)$. By \cite[Proposition 3.12]{LS},
$(S_1,\circ)$ is either trivial or has finite index in
$(S,\circ)$.
On the other, we claim that
\beq
\label{findex2}
\oo{uS_1} = \oo{vS_1} \iff u\pi = v\pi
\eeq
holds for all $u,v \in S$. The direct implication follows from $S_1\pi
= 1$. Conversely, assume that $u\pi = v\pi$. Then $\oo{v\inv u} \in
\oo{S_{v\inv\pi}S_{u\pi}} = S_1$ and so $u \in \oo{vS_1}$ and
$\oo{uS_1} \subseteq \oo{VS_1}$. By symmetry, we get $\oo{uS_1} =
\oo{VS_1}$ and so (\ref{findex2}) holds.

It is now straightforward to check that
$$\begin{array}{ccc}
(S,\circ)/(S_1,\circ)&\to&G\\
\oo{uS_1}&\mapsto&u\pi
\end{array}$$
is a group isomorphism.
Hence either $G \cong (S,\circ)$  is a free group of finite rank, or
$G \cong (S,\circ)/(S_1,\circ)$ is a finite group.

(ii) $\Rw$ (iii). Immediate from the proof of Proposition
\ref{freefin}.

(iii) $\Rw$ (i). Assume that $S = R_A$ is a Stallings section for the m-epi
$\pi:\wt{A}^* \to G$. Let $u \in S_g$ and $v \in S_h$ for some $g,h
\in G$. Since $\oo{uv}\pi = (uv)\pi = gh$, we get  $\oo{uv} \in
S_{gh}$ and so $\oo{S_gS_h} \subseteq S_{gh}$. Therefore
$S_{gh} = \oo{S_gS_h}$ and so $R_A$
satisfies (S2').
\qed


\section*{Acknowledgments}

The first author acknowledges support from the European Regional Development Fund through the programme COMPETE and by
the Portuguese Government through FCT (Funda\c c\~ao para a Ci\^encia e a Tecnologia) under the project
PEst-C/MAT/UI0144/2011. The second and third authors acknowledge support from the projects MTM2010-19938-C03-01 and
MTM2008-01550, respectively, of the Ministerio de Ciencia e Innovaci\'{o}n of the spanish government.

\end{document}